\documentclass[graybox,envcountsame,envcountresetsect]{svmult-fgs7}

\usepackage{type1cm}        % activate if the above 3 fonts are
\usepackage{makeidx}         % allows index generation
\usepackage{graphicx}        % standard LaTeX graphics tool
\usepackage{multicol}        % used for the two-column index
\usepackage[bottom]{footmisc}% places footnotes at page bottom

\usepackage{newtxtext}       % 
\usepackage[varvw]{newtxmath}       % selects Times Roman as basic font

\usepackage{amsmath}
\usepackage{amsthm}
\usepackage{verbatim}

\makeindex             % used for the subject index
\newtheorem{theo}{Theorem}[section]

\newtheorem{cor}[theo]{Corollary}

\newcommand{\R}{\mathbb{R}}

\renewcommand{\dh}{\dim_{\rm {H}}}
\newcommand{\db}{\dim_{\rm {B}}}
\renewcommand{\dlb}{{\underline\dim}_{\rm {B}}}
\newcommand{\dub}{{\overline\dim}_{\rm {B}}}
\newcommand{\dpk}{\dim_{\rm {P}}}
\newcommand{\proj}{ \mbox{\rm proj}}

\newcommand{\rn}{\mathbb{R}^n}

\newcommand{\toitself}{\mathbin{\scalebox{.85}{%
    \lefteqn{\scalebox{.5}{$\blacktriangleleft$}}\raisebox{0.34ex}{$\supset$}}}}

\begin{document}

\title*{Seventy Years of Fractal Projections}
\author{Kenneth J. Falconer\orcidID{0000-0001-8823-0406} }
\institute{Kenneth J. Falconer \at School of Mathematics and Statistics,  University of St Andrews, North Haugh, St Andrews, Fife KY16 8NS, UK, \email{kjf@st.andrews.ac.uk}
}
\maketitle

\abstract*{Seventy years ago, John Marstrand published a paper which, among other things, relates the Hausdorff dimension of a plane set to the dimensions of its orthogonal projections onto lines. For some time this paper attracted little attention, but over the past 40 years Marstrand's projection theorems have become the prototype for many results in fractal geometry with numerous variants and applications and they continue to motivate leading research.}

\section{Introduction}
\setcounter{equation}{0}
\setcounter{theo}{0}
\setcounter{figure}{0}
At the conference {\it Fractal Geometry and Stochastics V} held in Tabarz in 2014, I gave a survey talk entitled `Sixty years of fractal projections', a version of which, written with Jon Fraser and Xiong Jin \cite{FFJ}, appeared in in the conference proceedings \cite{BFZ}.  This marked the sixtieth anniversary of the publication of John Marstrand's 1954 paper \cite{Mar} `Some fundamental geometrical properties of plane sets of fractional dimensions'  in the Proceedings of the London Mathematical Society.  For a long time the paper attracted little attention, but since the 1980s, Marstrand's projection theorems have become the prototype for  many results in fractal geometry with numerous variants and applications. This area is now more intensively researched than ever,  drawing on modern techniques from ergodic theory, CP processes, Fourier transforms, discretisation of problems and additive combinatorics, together with a lot of ingenuity.

My talk `Seventy years of fractal projections' at {\it Fractal Geometry and Stochastics VII} held in 2024 in Chemnitz highlighted some of the very considerable progress in the area over the past 10 years. This account also surveys some of this more recent work and is a sequel to \cite{FFJ} where many earlier results and references may be found. There are many excellent articles and books which provide other substantial overviews of projection properties, including references \cite{Fra2,Mat3,Mat4,Mat5,Mat2,Mat6,Shm}. Also, volumes of conference proceedings, in particular of the German `Fractal Geometry and Stochastics' and the French `Fractals and Related Fields' meetings, contain many enlightening surveys on this and many other aspects of fractal geometry.

\subsection{General remarks}
Most of this article concerns orthogonal projections of sets in the plane onto straight lines or from 
$\mathbb{R}^n$ onto $m$-dimensional subspaces.
In the plane we let $L_\theta$ be the line making  angle $\theta \in [0,\pi)$ with the $x$-axis, and write $\proj_\theta: \mathbb{R}^2 \to L_\theta$ for orthogonal projection onto $L_\theta$, see Figure 1. By slight abuse of notation we will write ${\mathcal L} $ for Lebesgue measure on any line $L_\theta$, identified with $\mathbb{R}$ in the obvious way.

 \begin{figure}[h]
\begin{center}
\includegraphics[scale=0.4]{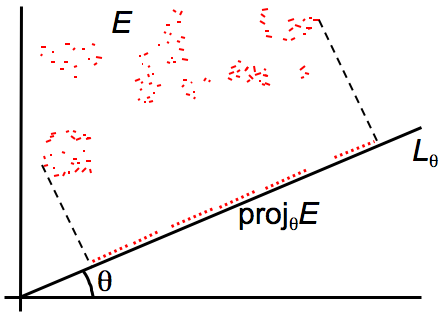}

Figure 1: Projection of a set $E$ onto a line in direction $\theta$.
\end{center}
\end{figure}
In the higher dimensional setting, we write $G(n,m)$ for the Grassmanian of $m$-dimensional subspaces of $\mathbb{R}^n$, where $1\leq m <n$. Then $G(n,m)$ is an $m(n-m)$-dimensional compact manifold which carries a natural invariant measure $\gamma_{n,m}$ that is locally equivalent to $m(n-m)$-dimensional Lebesgue measure. When we write `for almost all $V\in G(n,m)$' this is taken to be with respect to $\gamma_{n,m}$. For $V\in G(n,m)$ let $\proj_V : \mathbb{R}^n\to V$ be orthogonal projection onto  the $m$-dimensional subspace $V$, writing ${\mathcal L}^m$ for $m$-dimensional Lebesgue measure on $V$ which we identify with $\mathbb{R}^m$.

We make the convention that we use $\proj_\theta$ specifically for the planar setting and   $\proj_V$ for the more general $G(n,m)$ context.

Here our main interest is in projections of sets rather than measures, though many results for projections of sets have measure analogues, indeed many proofs concerning projections of sets depend on consideration of projections of suitable measures supported by the sets. Lower and upper Hausdorff and packing dimensions of measures can be defined in terms of local dimensions, and may all be equal for dynamically defined measures. Then often the dimension of projected measures can be expressed in terms of entropies and Lyapunov exponents. 

We will assume, sometimes without mentioning it specifically,  that the subsets of $\mathbb{R}^m$  with which we work  are all Borel or analytic, indeed for many results little is lost by assuming them to be compact. We avoid consideration of completely general sets that do not have analytic or measurable structure; they may behave strangely depending on the axioms of logic that are assumed. Nevertheless we note  that recently Lutz and Stull \cite{LS} have used effective dimension and computational complexity to obtain some results on projections with no requirement on sets to be Borel or analytic.

\subsection{Marstrand's projection theorem}
Marstrand's projection theorem in the plane \cite{Mar} may be stated as follows.
\begin{theo}\label{marthm}
Let $E \subset \mathbb{R}^2$ be a Borel or analytic set. Then

(i) $\dh \mbox{\rm proj}_\theta E \leq \min\{\dh E, 1\}$ with equality for almost all $\theta \in [0,\pi)$,
\smallskip

(ii) if $\dh E >1$ then ${\mathcal L} (\mbox{\rm proj}_\theta E) > 0$ for almost all $\theta \in [0,\pi)$.
\end{theo}

The  natural higher dimensional analogue of Marstrand's theorem was first presented by Mattila \cite{Mat} in 1975.

\begin{theo}\label{higherdim}
Let $E \subset \mathbb{R}^n$ be a Borel or analytic set. Then

(i) $\dh\mbox{\rm proj}_V E \leq \min\{\dh E, m\}$  with equality for $\gamma_{n,m}$-almost all $V \in G(n,m)$,
\smallskip

(ii) if $\dh E >m$ then ${\mathcal L}^m (\mbox{\rm proj}_V E) > 0$ for $\gamma_{n,m}$-almost all $V \in G(n,m)$. 
\end{theo}

Since orthogonal projection is a Lipschitz map that does not increase distances between points, the inequalities in (i) of Theorems \ref{marthm} and \ref{higherdim}   follow easily from the definition of Hausdorff measure and dimension, but showing that equality holds for almost all $\theta$ or $V$ is more involved. 
In his original paper, Marstrand  \cite{Mar} used intricate estimates involving plane geometry and measure theory, but in 1968  Kaufman \cite{Kau} gave a new proof of Theorem \ref{marthm}(i) using potential theory and of Theorem \ref{marthm}(ii)  using Fourier transforms before
Mattila \cite{Mat} used this approach in higher dimensions to obtain Theorem 1.2.

Briefly, the potential-theoretic method depends on the characterisation of Hausdorff dimension in terms of energy integrals. We write $\mathcal{M}(E)$ for the positive finite Borel measures supported on $E\subset \mathbb{R}^n$. For $s>0$ we write:
\begin{equation}\label{endef}
I^s(\mu) := \int \!\! \int \frac{d\mu(x)d\mu(y)}{|x-y|^s} 
\end{equation}
for the $s$-{\it energy} of the measure $\mu\in \mathcal{M}(E)$. Then
\begin{equation*}\label{endefhau}
\dh E = \sup\Big\{s :  \mbox{ there exists $\mu \in \mathcal{M}(E)$ such that $I^s(\mu) < \infty$}\Big\}.
\end{equation*}
Theorem \ref{marthm}(i) may then be proved by noting that for all $0<s< \dh E$ there is a measure $\mu\in \mathcal{M}(E)$ such that $I^s(\mu) <\infty$,   and verifying that  $\int I^s(\mu_\theta) d\theta \leq c_sI^s(\mu)$, where $\mu_\theta$ is the projection of the measure $\mu$ onto the line $L_\theta$ in direction $\theta$, given by $\mu_\theta(F) = \mu(\proj^{-1}_\theta F)$ for $F\subset L_\theta$. The higher dimensional case, Theorem \ref{higherdim}(i), may be proved in a similar way.

Theorem \ref{marthm}(ii) may be obtained using Fourier transforms. Writing $\widehat{\mu}(z):= \int e^{-2\pi i x.z} d\mu(x) \  (z\in \mathbb{R}^n)$ for the Fourier transform of $\mu$, applying Parseval's theorem and the convolution formulae to the energy expression in \eqref{endef} gives 
\begin{equation}\label{ftend}
I^s(\mu) = c_{n,s}\int_{\mathbb{R}^n} \frac{|\widehat{\mu}(z)|^2} {|z|^{n-s}}dz \qquad (0<s<n),
\end{equation}
where   $c_{n,s}$ depends only on $n$ and $s$, see \cite[Theorem 3.10]{Mat2}.
Then if $E\subset  \mathbb{R}^2$ and $\dh E >1$  we may find $\mu\in \mathcal{M}(E)$ and $1<s<\dh E$ such that $I^s(\mu) < \infty$. Writing \eqref{ftend} in radial coordinates this implies that $\int \!\! \int | \widehat{\mu_\theta}(t)|^2 dt\, d\theta <\infty$, so that for almost all $\theta$, the projected measure $\mu_\theta$ on $L_\theta$ is absolutely continuous with respect to $\mathcal{L}$ and is a square integrable function, so in particular $E$ has support of positive $\mathcal{L}$-measure.

Energy and Fourier methods and their generalisations and extensions have been used widely across fractal geometry and in particular in work on projections. Further examples will be encountered in this survey, but for thorough treatments see, for example, \cite{Fal, FFJ, Mat, Mat1,Mat4,Mat5,Mat2}

Marstrand's theorem is the prototype for a great deal of work in fractal geometry. If, in some setting, a dimensional property of parameterised images of a set is true for almost all parameters, perhaps for almost all $\theta$ or $V$, we say that a `Marstrand-type theorem holds'. We will see many instances of this in the sections that follow.

\section{Exceptional sets of projections}
\setcounter{equation}{0}
\setcounter{theo}{0}
\setcounter{figure}{0}

Whilst Marstrand's Theorem  gives that the Lebesgue measure  of the `exceptional set' of projection directions  $\theta$ for which $\dh \mbox{\rm proj}_\theta E < \min\{\dh E, 1\}$ has Lebesgue measure 0,  often the set of exceptional directions must be somewhat smaller than this. Indeed Kaufman's potential theoretic approach \cite{Kau} also yields an upper bound for the dimension of the exceptional set in Theorems \ref{marthm}(i) and \ref{higherdim}(i), namely that if $E \subset \mathbb{R}^n$  and $0\leq s< \dh E < m$ then
\begin{equation}\label{Kaut}
\dh\{V \in G(n,m) : \dh \mbox{\rm proj}_V E< s\} \leq m(n-m) -(m-s).
\end{equation}
(This bound is written in this form for comparison with $m(n-m)$, the dimension of the Grassmanian $G(n,m))$.
Using Fourier transforms, Falconer \cite{Fal1} in 1982 obtained  bounds for the exceptional sets of an alternative form: if $0\leq  s \leq m \leq \dh E$ then
\begin{equation}\label{Falt}
\dh\{V \in G(n,m) : \dh \mbox{\rm proj}_V E< s\} \leq \max\{m(n-m) -(\dh E-s), 0\},
\end{equation}
and if $m \leq \dh E \leq n$, 
\begin{equation}\label{Falt2}
\dh\{V \in G(n,m) :  \mathcal{L}^n(\mbox{\rm proj}_V E) =0\} \leq m(n-m) -(\dh E-m);
\end{equation}
the idea here is the greater the `excess dimension' $ \dh E - s$ or $ \dh E - m$, the smaller the set of exceptional $V$.
Estimates for the dimensions of exceptional sets in the spirit of \eqref{Kaut}, where there is a direct estimate, and \eqref{Falt}--\eqref{Falt2}, where the bound depends on by how much the dimension of a set exceeds a threshold,  have become known as `Kaufman-type' and `Falconer-type' estimates respectively.  

One further estimate, due to Peres and Schlag \cite{PSc} in 2000, completes this quartet.
If $\dh E > 2m$  then
\begin{equation}\label{int}
\dh\{V \in G(n,m) : \proj_V E \mbox{ \rm has  empty interior} \}\leq m(n-m)-  (\dh E-2m).\end{equation}

Examples suggested that the bounds for exceptional projections in \eqref{Kaut}-\eqref{Falt} were not in general optimal for all parameters, and various improvements were obtained particularly in the plane case, for example
Bourgain \cite{Bou2} used discretised methods to show  that when $ \dh E < 1$, 
\begin{equation}\label{Bour}
\dh\{\theta : \dh \mbox{proj}_\theta E< \textstyle{\frac{1}{2}}\dh E\} =0.
\end{equation}
Such results and various examples led to `Oberlin's Conjecture' \cite{Obe} proposed by Oberlin in 2012 for  the optimal upper bound: for $E\subset \R^2$ and $0 <t\leq \dh \leq 1$,
$$\dh\big\{\theta: \theta : \dh \proj_\theta E< \textstyle{\frac12}(t+\dh E)\big\} \leq t,$$
or equivalently 
\begin{equation*}\label{Ober1}
\dh\{\theta : \dh \proj_\theta E< s\} \leq \max\{2s - \dh E,0\},
\end{equation*}
for $ 0\leq s \leq \min\{\dh E,1\}$, see Figure $2$. 
 \begin{figure}[h]\label{Obe}
\begin{center}
\includegraphics[width=0.55\linewidth]{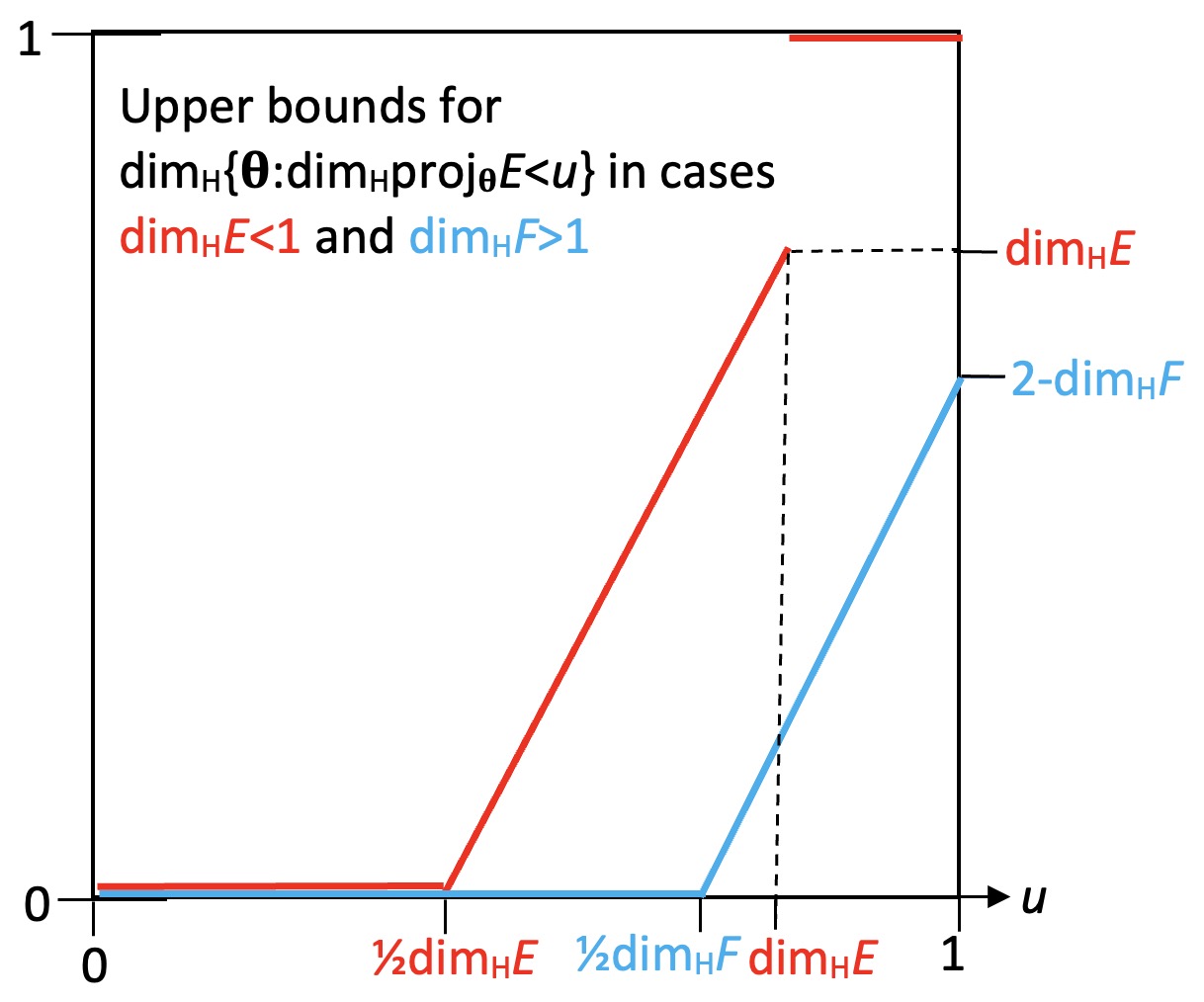}

Figure 2: Upper bounds for dimensions of exceptional directions of projections $\dh\{\theta:\dh\proj_\theta<s\}$ in cases $\dh E<1$ (broken line) and $\dh E>1$ (solid line). 
\end{center}
\end{figure}
Oberlin \cite{Obe} noted the analogy between this conjecture and the Furstenberg set conjecture. (This  concerns the minimum Hausdorff dimension of a set $E$ for which there is a family of lines $\mathbb{L}$  that all intersect $E$ in sets of Hausdorff dimension at least $s$, where  the line-set $\mathbb{L}$ has dimension at least $t$ in the natural parameterisation of lines.) Progress  towards Oberlin's conjecture was made by a number of authors including \cite{OO,Orp,OS,OS1} but it was finally fully established by Ren and Wang \cite{RW} in 2023 as a corollary to their proof of the Furstenberg set conjecture. 

\begin{theo}{\rm \cite{RW}}\label{obcon}
Let $E \subset \mathbb{R}^2$ be Borel and $ 0\leq s \leq \min\{\dh E,1\}$. Then 
\begin{equation}\label{Ober}
\dh\{\theta : \dh \proj_\theta E< s\} \leq \max\{2s - \dh E,0\},
\end{equation}
 and this bound is sharp.
\end{theo}

For projections of subsets of $\mathbb{R}^n$ onto $m$-dimensional subspaces for $n> 2$ knowledge is less complete. We have already noted  upper bounds for dimensions of exceptional projections in  \eqref{Kaut}--\eqref{int} but these bounds are not optimal for all parameters and better estimates are known in some cases, see \cite{Mat2}  for a discussion. 
For an example, He \cite{He} used a discretised approach to extend \eqref{Bour} to projections from $\mathbb{R}^n$ to $m$-dimensional subspaces:
$$\dh\big\{V\in G(n,m): \dh \mbox{proj}_V E\leq  \textstyle{\frac{m}{n}}\dh E\big\} \leq m(n-m)-1.$$

As an aside, we show that inequality \eqref{Ober} leads to a very short proof of the Erd\H{o}s-Volkmann conjecture, that there are no Borel subrings of $\mathbb{R}$, with the usual addition and multiplication,  that have Hausdorff dimension strictly between 0 and 1; this was proved in 2003 independently by Edgar and Millar \cite{EM} and Bourgain \cite{Bou}.  

For $a\in \mathbb{R}$ define  $\psi_a: \mathbb{R}^2 \to \mathbb{R}$  in coordinate form, also expressible as a `dot' product with the vector ${\bf \theta} =(a,1)$, by 
$$ \psi_a(x,y) = ax+y \equiv (x,y)\cdot \bf{ \theta}.$$
Geometrically, $ \psi_a$ is orthogonal projection onto a line in direction ${\bf \theta} =(a,1)$ together with a scaling of ratio $(a^2+1)^{1/2}$. 
Suppose   $R$ is a ring such that $0<\dh R<1$.  Let  $\epsilon>0$ be such that $\dh R +\epsilon <\min\{\dh (R\times R),1\}$.  (We may do this since $\dh (R\times R)\geq 2\dh R $, using the product rule for Hausdorff dimension that $\dh E + \dh F\leq \dh(E\times F)$ for Borel $E$ and $F$.)
Since  $R$ is a ring,  $\psi_a (R\times R) \subset R$ for all $a\in R$, so 
$$R\subset \{a: \dh\psi_a (R\times R) < \dh R +\epsilon\}.$$ 
Applying \eqref{Ober} with projections represented as $ \psi_a$, taking $E =R\times R$ and   $s= \dh R +\epsilon$,
\begin{align*}
\dh R&\leq \dh\{a: \dh\psi _a(R\times R) < \dh R+\epsilon  \}\\
&\leq\max\{2(\dh R +\epsilon) - \dh (R\times R),0\} \leq 2\epsilon,
\end{align*}
using the product rule again. This holds for arbitrarily small $\epsilon$, contradicting that $0<\dh R<1$.

\section{Other definitions of dimensions}\label{other}
\setcounter{equation}{0}
\setcounter{theo}{0}
\setcounter{figure}{0}

Whilst many problems involving projections were first considered for Hausdorff dimension, it is natural to ask similar questions using alternative definitions of dimension. In this section we look at projection properties relating to other forms of fractal dimension. 

\subsection{Box-counting dimensions}\label{box}

The {\it lower} and {\it upper box-counting dimensions} of  a non-empty and compact $E\subset  \mathbb{R}^n$ are given by 
$$
\dlb E\ =\ \liminf_{r\to 0} \frac{\log  N_r(E)}{-\log r} \ \text{ and }\ 
\dub E\ =\ \limsup_{r\to 0} \frac{\log  N_r(E)}{-\log r} 
$$
where $N_r(E)$ is the least number of sets of diameter $r$ covering $E$, see \cite{Fal,Mat1}. If $\dlb E =\dub E$ we write $\db E$  for the common value, termed the {\it box-counting}, {\it box}, {\it Minkowski} or  {\it Minkowski-Bouligand dimension}  of $E$. It is easy to see that $\dh E\leq \dlb E \leq \dub E$ for all $E\subset  \mathbb{R}^n$.

It is natural to ask whether there are Marstrand-type theorems for box-dimensions. J\"{a}rvenp\"{a}\"{a} \cite{Jar} constructed compact sets $E\subset \mathbb{R}^n$ with $\db E$ taking any prescribed value in $(0,n]$ and such that
$ \db \mbox{proj}_V E =  \db E\big/\big(1+ (1/m - 1/n)  \db E\big)$  for all $V\in G(n,m)$, 
which is strictly less than $\min\{ \db E, m\}$. Furthermore, in 1996 Falconer and Howroyd \cite{FH} showed that this lower bound was best possible and obtained a Marstrand-type theorem by showing that almost all projections of a set $E\subset \mathbb{R}^n$ must have the same lower, and the same upper, box dimensions. Their expressions for these values were awkward and difficult to work with and  Falconer revisited the question in 2020 \cite{Fal5, Fal4} using a potential-theoretic approach.

We define kernels  $\phi_r^s(x)$  for $s>0,\ 0<r<1$, $ x\in \mathbb{R}^n$ by
\begin{equation}\label{ker}
\phi_r^s(x)= \min\bigg\{ 1, \bigg(\frac{r}{|x|}\bigg)^{\! \! s}\bigg\}.
\end{equation}
The {\it capacity} $C_r^s(E) $ of a compact $E\subset \mathbb{R}^n$ with respect to  $\phi_r^s$  is given by
\begin{equation}\label{cap}
\frac{1}{C_r^s(E)}\  = \ \inf_{\mu \in {\mathcal M}_0(E)}\int\int \phi_r^s(x-y)d\mu(x)d\mu(y),
\end{equation}
where ${\mathcal M}_0(E)$ is the set of Borel probability measures supported by $E$. (For non-compact $E$ the capacity is taken to be the supremum of $C_r^s(F)$ over compact subsets $F$ of $E$.)
It may be shown that  for $E\subset \mathbb{R}^n$,
$$
c_1 C^s_r(E)\  \leq\ N_r(E)\  \leq \  
\left\{
\begin{array}{ll}
c_2\log(1 / r)\  C^s_r(E) & \mbox{ if } s=n    \\
c_2\ C^s_r(E)&  \mbox{ if } s>n   
\end{array}
\right. ,
$$
where $c_1, c_2$ depend on $s$ and the diameter of $E$.  In particular it follows that for $E\subset \mathbb{R}^n$,
\begin{equation*}\label{profilelow}
 \dlb^n E:= \dlb E =  \liminf_{r\to 0} \frac{\log  N_r(E)}{-\log r} =  \liminf_{r\to 0} \frac{\log  C_r^n(E)}{-\log r},
 \end{equation*}
 and
\begin{equation}\label{profile}
 \dub^n E:= \dub E =  \limsup_{r\to 0} \frac{\log  N_r(E)}{-\log r} =  \limsup_{r\to 0} \frac{\log  C_r^n(E)}{-\log r}.
 \end{equation}
 Taking $s=m$   in the kernel \eqref{ker} leads to the box dimensions of almost all projections.

\begin{theo}{\rm\cite{Fal5,Fal4}}\label{profilethm} Let $E \subset \mathbb{R}^n$ be non-empty and compact. 
Then for all $V \in G(n,m)$,
$$ \dlb \proj_V  E\  \leq \   \liminf_{r\to 0} \frac{\log  C_r^m(E)}{-\log r} =\dlb^mE,$$
and
$$ \dub \proj_V  E\  \leq \   \limsup_{r\to 0} \frac{\log  C_r^m(E)}{-\log r} =\dub^mE,$$
with equality for $\gamma_{n,m}$-almost all $V \in G(n,m)$. 
\end{theo}
Note that the almost sure dimensions of projections $\dlb^m E$ and  $\dub^m E$ may be strictly less than $\min\{\dlb E, m\}$ and  $\min\{\dub E, m\}$ respectively.

For $E \subset  \mathbb{R}^n$, when the capacities are defined with respect to the kernel  \eqref{ker}, we term
\begin{equation}\label{bprofdef}
\dlb^sE = \liminf_{r\to 0} \frac{\log  C_r^s(E)}{-\log r}\ \text{ and }\ 
\dub^sE = \limsup_{r\to 0} \frac{\log  C_r^s(E)}{-\log r} 
\end{equation}  
 the {\it  lower} and {\it  upper $s$-box-dimension profiles} of $E$, which should be thought of as the `box-dimension of $E$ when regarded from an $s$-dimensional viewpoint'. There is a parallel with Hausdorff dimensions, where one might define a `dimension profile' simply as 
$\dh^s E = \min\{s,\dh E\}$ which, by Marstrand's theorem, gives the Hausdorff dimension of projections onto  almost all $s$-dimensional subspaces.

The underlying reason why the kernels \eqref{ker} are central in studying projections is that there are numbers $a_{n,m}>0$ depending only on $n$ and $m$ such that for $x,y \in \rn, r>0$,
\begin{equation*}
\phi_r^m (x-y)\ \leq \ \gamma_{n,m}\big\{V\in G(n,m):  |\proj_Vx-\proj_Vy| \leq r\big\}
 \leq\ a_{n,m}\,\phi_r^m (x-y)
\end{equation*}
see \cite{Mat1,Mat2}.
This enables the integral over $V$ of the $m$-energies of the projected measures $\mu_V$ to be bounded
by the $m$-energy of $\mu$, for  suitable measures $\mu$ on $E$.

However, box dimensions cannot drop too much under projection: for almost all $V\in G(n,m)$
\begin{equation}\label{boxbounds} 
\frac{\dub E} {1+ (1/m - 1/n)  \dub E} \leq  \dub \mbox{proj}_V E \leq \min\{\dub E,m\}
\end{equation}
with similar sharp inequalities for $\dlb$, see \cite{Fal5, Fal4,FH,Jar}.  
 
As well as giving the almost sure box dimensions of the projections, the profiles provide upper bounds for the dimension of the exceptional sets of directions for which the dimensions fall below the almost sure value. In the next theorem the Kaufman-type bound (i) is analogous to \eqref{Kaut} for Hausdorff dimension, and the Falconer-type bound (ii) should be compared with \eqref{Falt}; these bounds are unlikely to be sharp in general.

\begin{theo}{\rm \cite{Fal5,Fal4}}\label{mainA}
Let $E\subset \mathbb{R}^n$ be a non-empty bounded Borel set. Then\\
(i) for $0\leq s\leq m$, 
$$\dh\{V \in  G(n,m): \dub \proj_V E \,  < \,   \dub^s E\}\   \leq\  m(n-m) - (m-s),$$
(ii) for $0\leq \gamma \leq n-m$,
$$\dh \{V \in  G(n,m) : \dub \proj_V E \,  < \,   \dub^{m+\gamma}E - \gamma \} \leq\  m(n-m) -\gamma.$$
(The  bound in (ii) is trivial unless  $0\leq \dub^{m+\gamma}E - \gamma\leq m$.) 
Directly analogous inequalities hold with $\dub E$ and $ \dub^s E$ replaced by $\dlb E$ and $ \dlb^s E$ respectively.
\end{theo}

The following more tractable but weaker bound is obtained in \cite{Fal5,FH2} using estimates of the kernels.

\begin{cor}\label{mainAcor}
Let $E\subset \mathbb{R}^n$ be a non-empty bounded Borel set. Then
 for $0< s\leq m$, 
$$\dh\Bigg\{ V\in G(n,m) : \dub \proj_V E <   \frac{\dub E} {1+ (1/s - 1/n)\dub E} \Bigg\} \leq m(n-m) - (m-s),$$
with a similar estimate where $\dub E$ is replaced by $\dlb E$.
\end{cor}

\subsection{Packing dimension}\label{pack}

Packing measures and packing dimension were  introduced  by Tricot \cite{Tri}  in 1982 as a sort of dual to their Hausdorff counterparts, see \cite{Fal,Mat1}. Whilst packing measures require an extra step in their definition, the gap of over sixty years between the two concepts seems with hindsight very surprising. As with Hausdorff dimension, packing dimension of a set $E$ is formally defined as the number $s$ at which the $s$-dimensional packing measure of $E$ changes from $\infty$ to 0.  We recall that $\dh E\leq \dpk E\leq \dub E$ for all $E\subset  \mathbb{R}^n$.

Following Ren and Wang's \cite{RW} proof of Oberlin's conjecture, Theorem \ref{obcon}, we immediately get a bound for projections in the plane with exceptional small packing dimensions: for $E\subset \mathbb{R}^2$,
$$\dh\{\theta: \dpk \proj_\theta E \leq s\} \leq \max\{2s - \dh E,0\}.$$
This supersedes earlier bounds in \cite{Orp7,OS}.

Many projection properties for box-counting dimensions transfer relatively easily to packing  dimensions since the packing dimension of a Borel or analytic $E\subset \mathbb{R}^n$ may alternatively be expressed in terms of the upper box dimensions of sets in countable coverings of $E$:  
$$\dpk E = \inf \Big\{ \sup_{1\leq i <\infty}\dub E_i : E \subset \bigcup_{i=1}^\infty E_i \text{ with $E_i$ compact}\Big\}.$$
In particular we may define the {\it packing dimension profile} of $E\subset \mathbb{R}^n$ for $s \geq 0$ by
$$\dpk^sE = \inf \Big\{ \sup_{1\leq i <\infty}\dub^s E_i : E \subset \bigcup_{i=1}^\infty E_i \text{ with $E_i$ compact}\Big\},$$
where  $\dub^s$ is as in \eqref{bprofdef}, so immediately  $\dpk E = \dpk^nE$.

With packing dimension defined in this way, the direct analogues of Theorem \ref{profilethm} and  the bounds \eqref{boxbounds} on dimensions of projections are easily established with  $\dub E$ and $ \dub^s E$ replaced by $\dpk E$ and $ \dpk^s E$ respectively, see \cite{Fal5, Fal4, FH2}.

\begin{theo}\label{packdim}
Let $E \subset \mathbb{R}^n$ be a Borel or analytic set. Then
$$\dpk \mbox{\rm proj}_V E \leq \dpk^m E,$$
with equality for $\gamma_{n,m}$-almost all $V \in G(n,m)$. 
\end{theo}

As with box-dimensions, the packing dimension profiles provide upper bounds for the Hausdorff dimension of the exceptional set of directions for which the packing dimension falls below the almost sure value. Again bounds on the dimensions of exceptional directions are given by Theorem \ref{mainA}(i) and Corollary \ref{mainAcor}  where $\dub E$ and $ \dub^s E$ are replaced by $\dpk E$ and $ \dpk^s E$.

Since their introduction, packing dimension profiles have cropped up in other contexts, notably to give the almost sure packing dimension of images of sets under fractional Brownian motion \cite{KX,Xia}.

In the case of measures as opposed to sets, there is a refined lower bound for the packing dimension of a measure for projections in almost all directions that incorporates both the Hausdorff and packing dimensions of the measure, see \cite{FM}.

Orponen \cite{Orp} notes that if  $\theta_1\neq \theta_2$ then the product rule for packing dimensions implies that  $\dpk E \leq \dpk \proj_{\theta_1}E  + \dpk \proj_{\theta_2}E$ for any Borel $E\subset \mathbb{R}^2$, from which it is immediate that, if $\dpk E>0$,
\begin{equation}\label{packhalf}
\dpk \proj_\theta E < \textstyle{\frac12} \dpk E
\end{equation} 
for at most one direction $\theta$. 
However, replacing `$<$'  by `$\leq$' he gives an example of a compact  $E\subset \mathbb{R}^2$ such that
$$\dpk\{\theta: \dpk \proj_\theta E \leq \textstyle{\frac12} \dpk E\} = 1$$ 

Since the packing dimension of a set is always at least its Hausdorff dimension, obtaining bounds for the packing dimension of exceptional sets of directions is harder than for  Hausdorff dimension bounds.  
Orponen \cite{Orp} has some nice bounds in the plane: for analytic $E\subset \mathbb{R}^2$ with $\dh E\leq 1$,
$$
\dpk\{\theta: \dpk \proj_\theta E\leq s\} \leq \  
\left\{
\begin{array}{ll}
\frac{s \dh E}{\dh E + s(\dh E -1)} & \mbox{ if } 0 \leq s\leq \dh E,    \\
\\
\frac{(2s - \dh E)(1-\dh E)}{\dh \!E /2} +s & \mbox{ if }  \frac12 \!\dh E \leq s\leq \dh E   
\end{array}
\right. .
$$

\subsection{Intermediate dimensions}\label{inter}

Intermediate dimensions were introduced by Falconer, Fraser and Kempton \cite{FFK} in 2020 to interpolate between Hausdorff dimensions and box dimensions. For $0\leq \vartheta\leq 1$ we define the {\it upper $\vartheta$-intermediate dimension} of a non-empty bounded set $E\subset \mathbb{R}^n$ by
\begin{align*}
\overline{\dim}_{\vartheta}E  =  \inf &\big\{ s\geq 0  :  \mbox{\rm for all $\epsilon >0$ and all  sufficiently small $\delta>0$}\\
 & \mbox{ \rm there is a cover $ \{U_i\} $ of $E$ s.t. $\delta^{1/\vartheta} \leq  |U_i| \leq \delta$ and $\sum |U_i|^s \leq \epsilon$}  \big\}.
\end{align*}
The {\it lower  $\vartheta$-intermediate dimension }$\underline{\dim}_{\vartheta}E$ is defined in the same way except the conditions are only required to hold for  a sequence of $\delta$ approaching $0$. Below we just consider the upper definition, but everything goes across directly to the lower intermediate dimension case.

The intermediate dimension $\overline{\dim}_{\vartheta}E$  {\it interpolates} between Hausdorff and box dimensions, that is  $\vartheta\mapsto \overline{\dim}_{\vartheta}E$ is increasing for $\vartheta \in [0,1]$ and
$$\dh E = \overline{\dim}_0 E \leq \overline{\dim}_{\vartheta}E
\leq \overline{\dim}_{1}E  = \dub E.$$ 
The function $\vartheta\mapsto \overline{\dim}_{\vartheta}E$ is continuous on $(0,1]$ and may or may not be continuous at 0.
The basic properties of intermediate dimensions are described in \cite{Fal6,FFK} and for $\vartheta \in (0,1]$ behave more like box-dimensions than Hausdorff dimension.

Analogously to box-counting dimension \eqref{bprofdef}, for $E \subset \mathbb{R}^n$, $s\geq 0$ and $\vartheta \in [0,1]$, we define the dimension profiles
$\overline{\dim}_{\vartheta}^s E$
in terms of energy integrals with respect to the kernels 
$${\phi}_{r, \vartheta}^{s, m}(x)
 = \begin{cases} 
      1 & 0\leq |x| < r \\
      \big(\frac{r}{|x|}\big)^s & r\leq |x| < r^\vartheta   \\
      \frac{r^{\vartheta(m-s) + s}}{|x|^m}\ & r^\vartheta \leq |x|.
   \end{cases} 
$$
Then for $E \subset \mathbb{R}^n$ and $\vartheta \in [0,1]$,
$\overline{\dim}_{\vartheta} E = \overline{\dim}_{\vartheta}^n E$.
Moreover, as with box-dimensions, other values of $s$ give the intermediate dimensions of projections of $E$. The following Marstrand-type theorem was obtained by Burrell, Falconer and Fraser \cite{BFF}, see also \cite{Fal6, Fra2}.

\begin{theo}\label{intproj}
Let $E \subset \mathbb{R}^n$ be a Borel or analytic set. Then for $\gamma_{n,m}$-almost all $V \in G(n,m)$,
$$ \overline{\dim}_{\vartheta} \proj_V E =\overline{\dim}_{\vartheta}^m E.
$$
for all $\vartheta \in [0,1]$.
\end{theo}

\subsection{Assouad dimension}\label{ass}
Recall that we write  $B(x,r)\subset \mathbb{R}^n$ for the ball of centre $x$ and radius $r$ and $N_r(F)$ for the least number of sets of diameter $r>0$ that can cover $F \subset \mathbb{R}^n$. 
The {\it Assouad dimension} of $E \subset \mathbb{R}^n$ is given by
\begin{align*}
\dim_{\rm A}E = \inf \bigg\{s :& \text{ there exists $c>0$ s.t.  for all }0<r<R \\
&\text{ and } x\in E, \ N_r(B(x,R)\cap E) \leq c \bigg(\frac{R}{r}\bigg)^s\bigg\}.
\end{align*}
Assouad dimension has attracted a great deal of interest in recent years; its character differs from most other fractal dimensions in that it has a local nature, highlighting the `thicker' parts of a set. Assouad dimension and its variants have interesting relationships with other forms of dimension, and are useful tools in embedding theory and in studying the regularity of mappings. The book by Fraser \cite{FraBk} provides a comprehensive treatment.

Assouad dimension does not share all the `standard' properties of other fractal dimensions. In particular the Assouad dimension of a Lipschitz image of a set can be less than or greater than that of the set, although Assouad dimension is preserved under bi-Lipschitz mappings. The absence of a Lipschitz property is manifest in its behaviour under projections.

There is no `almost sure' Marstrand-type result for Assouad dimension.   Let $E_s$ be the  $s$-dimensional right Sierpi\'{n}ski triangle, that is the self-similar set based on three homotheties (similarities without rotations) of ratio $3^{-1/s}$, where $0<s<1$, so $\dh E_s = {\dim}_{\rm A}E_s =s$, see Figure 3.  
\begin{figure}[h]
\begin{center}
\includegraphics[width=0.55\linewidth]{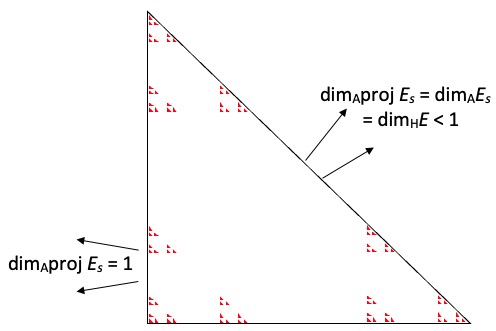}

Figure 3: The set $E_s$ with projections of Assouad dimension 1 or $s<1$  in different ranges of directions.
\end{center}
\end{figure}

\begin{theo}\label{assd} 
Let $\log 3/\log 5 <s<1$. Then there is $\epsilon>0$ such that
${\dim}_{\rm A}\proj_\theta E_s=1$ for almost all  $\theta \in (-\epsilon, \epsilon)$ and ${\dim}_{\rm A}\proj_\theta E_s=s<1$ for all $\theta \in (\pi/4 -\epsilon, \pi/4+ \epsilon)$.
\end{theo}
\noindent This example, due to Fraser and Orponen \cite{FO}, depends on the fact that a self-similar $E\subset \mathbb{R}^2$ defined by contracting homotheties   has
$$
{\dim}_{\rm A}\mbox{proj}_\theta E =
\left\{
\begin{array}{ll}
 \dh \proj_\theta E& \text{ if  } \mathcal{H}^{\dh {\rm proj}_\theta E}(\proj_\theta E) >0 \\
 1 & \text{ if  }\mathcal{H}^{\dh {\rm proj}_\theta E}(\proj_\theta E)=0   
\end{array}
\right.,
$$
where $\mathcal{H}^s$ is $s$-dimensional Hausdorff measure, so the construction is equivalent to finding a set of homotheties for which the IFS attractor $E$ has  projections with positive or zero ($\proj_\theta E$)-dimensional Hausdorff measure in the appropriate directions. 

Fraser and K\"{a}enm\"{a}ki\cite{FK} showed that Assouad dimensions of projections can be even more bizarre: for any upper semi-continuous $g :[0,\pi)\to [0,1]$, there exists a compact $E\subset \mathbb{R}^2$ such that 
${\dim}_{\rm A}\mbox{proj}_\theta E = g(\theta)$ for all $\theta$.

At least, in almost all directions, Assouad dimension does not drop under projection. Fraser and  Orponen \cite{FO} and Fraser  \cite{Fra1} showed that for $E\subset \mathbb{R}^n$ 
\begin{equation}\label{assproj}
{\dim}_{\rm A} \proj_V E \geq \min\{m, {\dim}_{\rm A} E\} 
\end{equation}
for almost all  $V \in G(n,m)$ and Orponen \cite{Orp6} obtained a much stronger bound on the exceptional set of directions in the plane: 
$$\dh\big\{\theta : {\dim}_{\rm  A}\mbox{proj}_\theta E< \min\{1, {\dim}_{\rm A} E\}\big\}=0. $$ 
Just as for packing dimension \eqref{packhalf}, the inequality for the Assouad dimension of a product implies that for $E\subset \mathbb{R}^2$ that
${\dim}_{\rm A} \proj_\theta E < \textstyle{\frac12} {\dim}_{\rm A} E$ can occur for at most one direction $\theta$.

Variants on Assouad dimension include  the {\it Assouad spectrum } $\dim_{\rm A}^\vartheta E\ (0<\vartheta<1)$ of $E \subset \mathbb{R}^n$ given by
\begin{align*}\dim_{\rm A}^\vartheta E = \inf \bigg\{s : &\text{ there exists $c>0$ s.t.  for all }0<r<1\\
&\text{ and } x\in E, \ N_r(B(x,r^\vartheta)\cap F) \leq c \bigg(\frac{r^\vartheta}{r}\bigg)^s\bigg\},
\end{align*}
and the related {\it  quasi-Assouad dimension } which may be expressed as
$$\dim_{\rm qA} E = \lim_{\vartheta\to 1} \dim_{\rm A}^\vartheta E.$$
 It seems unknown whether there are Marstrand-type projection results for $\dim_{\rm  A}^\vartheta E$ for each $0< \vartheta<1$ and  for $\dim_{\rm qA} E$.

\subsection{Fourier dimension}\label{Four}

We define the {\it Fourier transform} of a finite measure $\mu$ on $\mathbb{R}^n$ by
$$\widehat{\mu}(z) = \int_{x\in \mathbb{R}^n} e^{-2\pi i x\cdot z} d\mu(x) \quad (z\in \mathbb{R}^n).$$
Then the {\it Fourier dimension } $\dim_{\rm F}E$ of $E \subset \mathbb{R}^n$ is given by
\begin{align*}
\dim_{\rm F}E= \sup \big\{&s\leq n: \text{ there exists }  c>0  \text{  and  }\mu\in \mathcal{M}(E)  \\
&\text{such that  } |\widehat{\mu}(z) | \leq c |z|^{-s/2} \text{ for all } z\in \mathbb{R}^n\big\},
\end{align*}
reflecting the rate of decay of the Fourier transform of measures supported by $E$. Then $\dim_{\rm F} E \leq \dh E$ by the Fourier characterisation \eqref{ftend} of $\dh E$, since for all probability measures $\mu$ and $0<t< s<\dim_{\rm F} E$,
$$  
 \int \frac{|\widehat{\mu}(z) |^2 dz}{|z|^{n-t}} \leq \int_{|z|<1} \frac{dz}{|z|^{n-t}} + \int_{|z|\geq 1} \frac{c^2 |z |^{-s} dz}{|z|^{n-t}} < \infty,$$
implying that $\dh E\geq t$ for all $t<\dim_{\rm F} E$. 

Fourier transforms behave well with respect to projections: for a measure $\mu$ on $ \mathbb{R}^n$
\begin{equation}\label{Fourproj}
\widehat{\mu_V}(z)= \widehat{\mu} (z) \ \text{ for } z \in V \in G(n,m),
\end{equation} 
where $\mu_V$ is the projection of the measure $\mu$ onto $V$ given by $\mu_V (F) = \mu(\proj_V^{-1}(F))$ for $F \subset V$. 
Then for every $V \in G(n,m)$,
\begin{equation}\label{Fourierin}
\min\{m, \dim_{\rm F} E\} \leq \dim_{\rm F}  \proj_V E \leq \dh  \proj_V E\leq \min\{m, \dh E\},
\end{equation}
where the left-hand inequality follows on considering decay rates in \eqref{Fourproj} and the right-hand one is the upper bound from Marstrand's  theorem.
In particular, if $E$ is a {\it Salem set}, that is if  $\dim_{\rm F} E = \dh E$, there is equality throughout \eqref{Fourierin}, so that 
$$\qquad\dim_{\rm F}  \proj_V E =\dh  \proj_V E\ = \min\{m, \dh E\} \text { for   all }  V\in G(n,m),$$
that is there are no exceptional directions.

In general, only a limited amount of information on projections can be gleaned from the Fourier dimension of a set alone, indeed 
it is unknown whether there is a Marstrand-type result for $\dim_{\rm F}\proj_V E$. Fourier dimension itself  does not give extra information about the exceptional directions for Hausdorff dimensions of projections. Fraser and de Orellana \cite{FdO} constructed compact sets $E  \subset \mathbb{R}^2$ with  $\dh E=   s\in (0,1]$ and $ \dim_{\rm F}E=t\in (s/2,s)$ such that 
$$
\begin{array}{rlll}
\{\theta: \dh \proj_\theta E \leq u\} &=& \varnothing & \text{ if  } u<t  \\
\dh \{\theta: \dh \proj_\theta E \leq u\} &=&2t-s & \text{ if  }u\geq t   
\end{array}
$$
so that the dimension of the exceptional set has a jump discontinuity at $t$ from 0 to the largest value permissible by \eqref{Ober}.

To provide more information,  
Fraser \cite{Fra3} introduced the {\it Fourier spectrum} $\dim_{\rm F}^\vartheta$ of a set $E \subset\mathbb{R}^n$ that interpolates between $\dim_{\rm F}E$ when $\vartheta =0$ and $\dh E$ when $\vartheta =1$, which leads to many further inequalities, see \cite{Fra3, Fra2, FO}. Bounds for the size of exceptional sets for the Fourier spectrum of projections may be obtained for $\vartheta \in (0,1]$: for all $0\leq s\leq \min\{m, s- \dim_{\rm F}^\vartheta E\}$:
$$\dh \{V\in G(n,m): \dim_{\rm F}^\vartheta \proj_V E \leq s\} 
\leq \max\bigg\{0, m(n-m) + \inf_{\vartheta'\in (0,1]}\frac{s- \dim_{\rm F}^{\vartheta'} E}{\vartheta'}\bigg\}.
$$                                                
Setting   $\vartheta=1$  gives  an upper bound for                                                                                                                         
$\dh \{V\in G(n,m): \dh \proj_V E \leq s\} $
which, by taking a suitable value of $\vartheta'$ in the infimum, can under some circumstances give a better bound than \eqref{Ober}, see \cite{FO}.

\section{Projections in restricted directions}
\setcounter{equation}{0}
\setcounter{theo}{0}
\setcounter{figure}{0}
A general question that has been around for many years asks
under what circumstances can one get projection results for projections onto families of lines or subspaces that form proper subsets of $G(n,m)$.  For instance, if
$\{V(t): t \in P\}$ is a smooth curve or submanifold 
 of $G(n,m)$ smoothly  parameterized by a set  $P \subset \mathbb{R}^k$, then what can we conclude about
 $\dh \proj_{V(t)} E$ for ${\mathcal L}^k$-almost all $t \in P$, where  ${\mathcal L}^k$ is $k$-dimensional Lebsegue measure?
 
 For an easy example, let $\theta: [0,1] \to S^2$ (the 2-sphere embedded in $\mathbb{R}^3$) be a smoothly parameterised curve of directions and write $\proj_{\theta(t)}$ for orthogonal projection onto the line in direction $\theta(t)$.  Then for almost all $0 \leq t \leq 1$ we have $\dh \proj_{\theta(t)} E\geq \min\{\dh E -1,0\}$, with $ {\mathcal L}^1(\proj_{\theta (t)} E)>0$  if $\dh E > 2$, since by
 \eqref{Falt} and  \eqref{Falt2}  the set of directions for which these lower bounds fail has Hausdorff dimension less than 1, the Hausdorff dimension of the curve.

The following lower bounds  were obtained J\"{a}rvenp\"{a}\"{a}, J\"{a}rvenp\"{a}\"{a} and Keleti \cite{JJK} for parameterized families of projections from  $\mathbb{R}^n$ to $m$-dimensional subspaces, see also \cite{JJLL}. 
  For $0<k<m(n-m)$ define the integers
$$p(l) = n-m -\bigg\lfloor \frac{k-l(n-m)}{m-l}\bigg\rfloor \qquad (l = 0,1,\ldots, m-1),$$
where the `floor' symbol `$\lfloor x\rfloor$' denotes the largest integer no greater than $x$.

\begin{theo}\label{projgeneral}
Let $P \subset \mathbb{R}^k$ be an open parameter set and let $E \subset \mathbb{R}^n$ be a Borel or analytic set. Let $\{V(t)\subset G(n,m): t \in P\}$ be a  family of subspaces such that $V$ is $C^1$ with the derivative $D_t V(t)$ injective for all $t \in P$. Then, for all $l=0,1,\ldots,m$ and ${\mathcal L}^k$-almost all $t \in P$,
$$
\dh\proj_{V(t)} E \geq  
\left\{
\begin{array}{llrl}
 \dh E-p(l) &\ \mbox{if} &  p(l) + l   & \leq \dh E \leq p(l) + l +1  \\
l+1 &\ \mbox{if}  & \ p(l) + l +1  &\leq \dh E \leq p(l+1) + l +1
\end{array}
\right. .
$$
Moreover, if $\dh E > p(m-1) + m$ then  ${\mathcal L}^m ( \proj_{V(t)} E) > 0$ for ${\mathcal L}^k$-almost all $t \in P$.
\end{theo}

These are the best possible bounds for general parameterised families of projections.
The same paper \cite{JJK} includes extensions of these results to smoothly parameterised families of $C^2$-mappings.

Better bounds may be obtained  if there is curvature in the mapping $t \mapsto V(t)$. This is a difficult area but there has been impressive progress in recent years on projections from $\mathbb{R}^3$ to lines or planes, with many contributions including \cite{FO1,Fra2,GGM,Har2,Har, Har1,He,KOV,Obe, Obe1, OO,Orp,Orp1,OS, OS1,OV,PYZ} bringing in a host of new methods, including techniques from Fourier and harmonic analysis, discretisation, additive combinatorics and decoupling. 

Let  $\theta: [0,1] \to S^2$ be a family of directions given by  a $C^3$-function $\theta$. We say that the curve of directions  is {\it non-degenerate} if 
$$\mbox{span }\{\theta(t), \theta'(t), \theta''(t)\} = \mathbb{R}^3 \ \mbox{ for all  } t \in [0,1],$$
for instance the curve $t\mapsto \frac{1}{\sqrt{2}}(\cos t, \sin t, 1)$ satisfies this condition. As before $\proj_{\theta(t)}$ denotes orthogonal projection onto the line in direction $\theta(t)$.
The following theorem, which includes a   Marstrand-type theorem and bounds on the dimension of $t$ for which images under $\proj_{\theta(t)}$ are  exceptionally small, summarises recent progress. 

\begin{theo}\label{projline}
Let $E \subset \mathbb{R}^3$ be a Borel or analytic set and let $\theta: [0,1] \to S^2$ be a non-degenerate family of directions. Then 
\begin{equation}\label{projlineeq}
\dh\proj_{\theta(t)} E = \min\{\dh E, 1\}\quad \text{ for almost all } t\in [0,1],
\end{equation}
and if $\dh E >1$ then
\begin{equation}\label{resmes}
{\mathcal L} (\proj_{\theta(t)} E)>0 \quad \text{ for almost all } t\in [0,1].
\end{equation}
Concerning exceptional directions, if $0\leq s\leq \min\{\dh E, 1\}$ then
\begin{equation}\label{PYZ}
 \dh\{t : \dh \proj_{\theta(t)} E <s\} \leq s,
\end{equation}
and if $0\leq s\leq 1$ then 
\begin{equation}\label{GGM}
 \dh \{t: \proj_{\theta(t)} E <s\} \leq \max\big\{0, 1 + {\textstyle \frac12} (s-\dh E) \big\},
\end{equation}
and
\begin{equation}\label{resmes1}
\dh\{t: {\mathcal L} (\proj_{\theta(t)} E) = 0 \} \leq  \textstyle{\frac13}(4-\dh E).
\end{equation}
\end{theo}

\noindent Pramanik, Yang and Zahl \cite{PYZ} proved \eqref{PYZ} using a circular maximal function and Gan, Guth and Maldague \cite{GGM} used a decoupling method to obtain  \eqref{GGM}. The Marstrand-type result \eqref{projlineeq} follows from either of these.  
Harris \cite{Har, Har1} separately obtained  \eqref{resmes} and \eqref{resmes1}.

As a complementary operation we may project from $\mathbb{R}^3$ onto planes rather than lines and again there has been considerable recent progress on the almost sure dimensions of projections onto parameterised families of planes, alongside bounds for exceptional sets. In this case we write $\proj_{V_{\theta}}: \mathbb{R}^3\to V$ for projection onto the plane perpendicular to the vector $\theta \in S^2$.
Work including \cite{FO1,Orp1,OV} led to the following result proved by Gan et al \cite{GGGHMW} using the `high-low' method. The values \eqref{projplaneeq}--\eqref{GGG}  should be compared with the non-restricted projection results in Theorem \ref{higherdim} and \eqref{Kaut}--\eqref{Falt2}.

\begin{theo}\label{projplane}
Let $E \subset \mathbb{R}^3$ be a Borel or analytic set and let $\theta: [0,1] \to S^2$ be a non-degenerate family of directions. 
Then 
\begin{equation}\label{projplaneeq}
\dh\proj_{V_{\theta(t)}} E = \min\{\dh E, 2\} \quad \text{ for  almost all } t\in [0,1].
\end{equation}
Moreover, if $\dh E >2$ then
\begin{equation}\label{projlplanex}
{\mathcal L^2} (\proj_{\theta(t)} E)>0 \quad \text{ for almost all } t\in [0,1].
\end{equation}
For  exceptional directions, if $0\leq s\leq 1$ then 
\begin{equation}\label{GGG}
 \dh\{t : \dh \proj_{\theta(t)} E <s\} \leq \max\{ 1 + s-\dh E, 0 \}.
\end{equation}
\end{theo}

A higher dimensional version of this theorem was recently developed by some of these authors \cite{GGW}. They showed that if $\theta: [0,1] \to \mathbb{R}^n$ is a smooth non-degenerate curve and $m\leq n$, then the projection of $E$ to the $m$th order tangent space of $\theta$ at $t\in [0,1]$ has Hausdorff dimension $\min\{\dh E, m\}$ for almost all $t$. This turns out to have applications to Diophantine approximation. The Oppenheim conjecture, formulated by Alex Oppenheim in 1929 and originally proved by Margulis \cite{Marg} in 1987, states that if an indefinite nondegenerate quadratic form of three or more variables is not proportional to a rational quadratic form, then its set of values at the integers is dense in the real numbers.  A deep paper \cite{LMWY} obtains a quantitative version of the Oppenheim conjecture with a polynomial error term, improving earlier estimates;  a key step in the proof uses this restricted projection result taking  $n=5$ and $m=2$.

Also in higher dimensions, Mattila \cite{Mat7} has obtained some nice inequalities for a restricted class of projections from $\mathbb{R}^{2n}$ to $\mathbb{R}^{n}$ for $n\in \mathbb{N}$.

\section{Projections of IFS attractors: self-similar and self-affine sets}\label{secss}
\setcounter{equation}{0}
\setcounter{theo}{0}
\setcounter{figure}{0}

One of the drawbacks of the projection theorems and exceptional set results  is that they tell us nothing about the dimension or measure of the projection of a set $E$ in any prescribed direction. However, if  $E$ has some regularity, for example  some form of self-similarity or self-affinity, then more can often be said. 
There has been considerable recent interest in examining the dimensions of projections in specific directions for particular sets or classes of sets, and especially in finding sets for which the conclusions of Marstrand's theorems are valid for all, or virtually all, directions. 

Recall that an {\it iterated function system} (IFS)  is a family of  contractions $\{f_1, \ldots,f_k\}$ with $f_i : \mathbb{R}^n\to \mathbb{R}^n$ (or sometimes with $f_i : D\to D$ for some closed domain $D\subset  \mathbb{R}^n$). An IFS determines a unique non-empty compact  $E\subset \mathbb{R}^d$ (or $E\subset D$) such that 
\begin{equation}\label{attractor}
E= \bigcup_{i=1}^k f_i(E),
\end{equation}
called the {\it attractor} of the IFS, see, for example,  \cite{Fal,Hut}. 
Moreover, $E$ can be realised by a hierarchical construction. Let $A\subset \mathbb{R}^n$ be nonempty and compact such that $f_i(A)\subset A$ for all $1\leq i\leq k$. Writing $A_{i_1,i_2,\ldots,i_j}= f_{i_1}\circ f_{i_2}\circ \cdots \circ f_{i_j}(A)$ and $F_j = \bigcup_{1\leq i_1,\ldots, i_j \leq k} A_{i_1,i_2,\ldots,i_j}$ the attractor is given by 
\begin{equation}\label{attint}
E= \bigcap_{j=0}^\infty F_j.
\end{equation}

In this section we consider several important classes of attractors, see Figure 4,  including self-similar and self-affine attractors for which the recent book by B\'{a}r\'{a}ny, Simon and Solomyak \cite{BSS} is excellent reference.

\begin{figure}[h]\label{IFS}
\begin{center}
\includegraphics[width=0.90\linewidth]{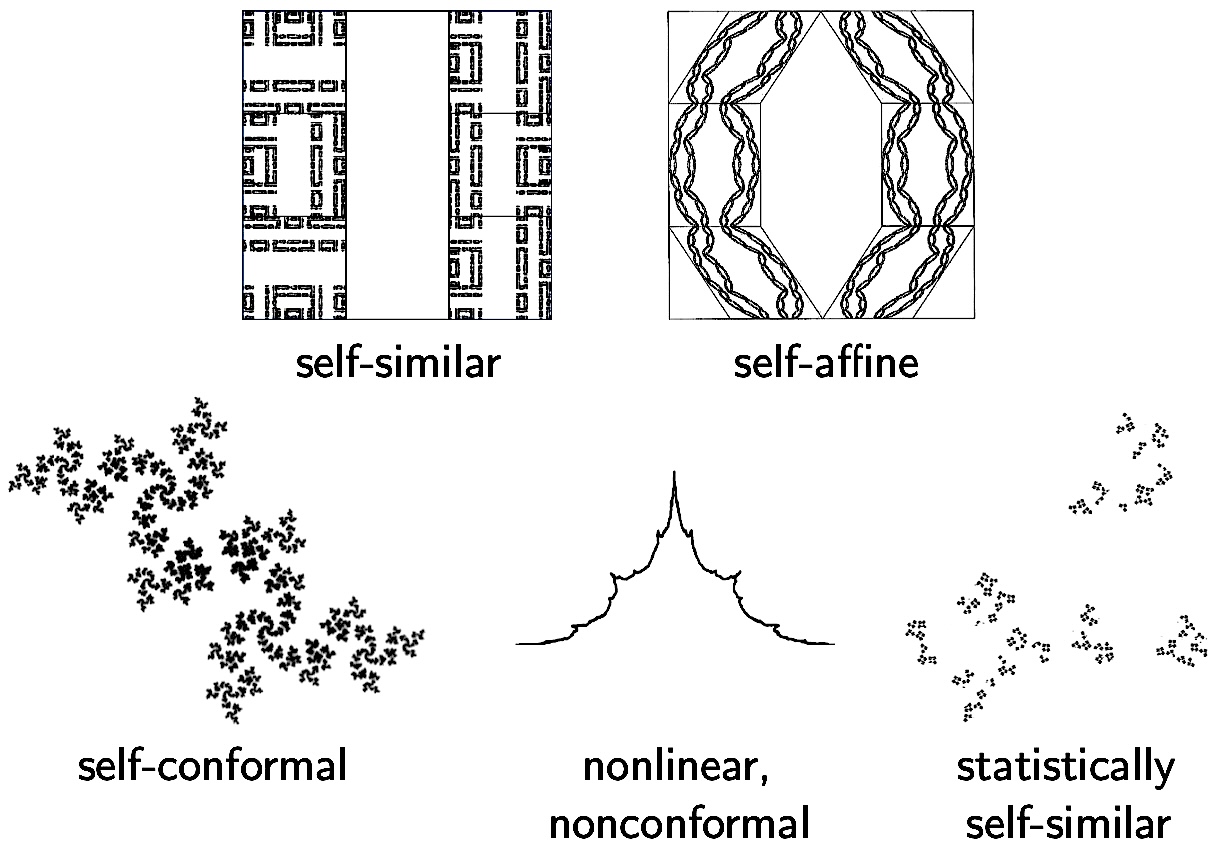}
\medskip

Figure 4: Attractors of IFSs of mappings of different types.
\end{center}
\end{figure}

\subsection{Self-similar sets}\label{sss}

If the $f_i$ in the IFS are all similarities, that is of the form
\begin{equation}\label{simform}
f_i(x) = r_i O_i (x) + a_i \quad (1\leq i \leq k),
\end{equation}
where $0<r_i<1$ is the contraction ratio, $O_i$ is an orthonormal map, i.e. a rotation or reflection, and $a_i$ is a translation,
the attractor $E$ is termed {\it self-similar}. 

An IFS of similarities (or sometimes  more general mappings) satisfies the {\it strong separation  condition} (SSC) if the union (\ref{attractor}) is disjoint,  the {\it open set condition} (OSC) if there is a non-empty open set $U$ such that $\bigcup_{i=1}^k f_i(U)\subset U$ with this union disjoint, and the {\it strong open set condition} (SOSC) if such a $U$ can be chosen so $U\cap E \neq \varnothing$. 

If either SSC or OSC hold then $\dh E = s$, where $s$ is the {\it similarity dimension} given by $\sum_{i=1}^k r_i^s=1$, where $r_i$ is the similarity ratio of $f_i$, and moreover $0<{\mathcal H}^s(E) < \infty$ where ${\mathcal H}^s$ is $s$-dimensional Hausdorff measure. The {\it rotation group} $G= \langle O_1, \ldots, O_k \rangle$ generated by the orthonormal components of the similarities plays a crucial role in the behaviour of the projections of self-similar sets, with conclusions depending on whether  $G$ is finite or not.

It is easy to construct self-similar sets with a finite rotation group $G$ for which the conclusions of Marstrand's theorem fail in certain directions. For example, let 
$f_1,\ldots,f_4$ be  homotheties (that is similarities with $O_i$ the identity in \eqref{simform}) of ratio $0<r<\frac{1}{4}$  that map the unit square $S$ into itself, each $f_i$ fixing one of the four corners. Then $\dh E = -\log 4 /\log r$,   but the projections of $E$  onto the sides of the square have  dimension $ -\log 2 /\log r$
 and onto the diagonals of $S$ have dimension $ -\log3 /\log r$, a consequence of the alignment of  the component squares $f_i(S)$ under projection. 
 
 In fact, when the rotation group is finite, there are always some projections under which there is a drop in dimension, as the following theorem of Farkas \cite{Far} asserts.
\begin{theo}\label{mesproj}
If $E\subset \mathbb{R}^n$ is self-similar with finite rotation group $G$ and similarity dimension $s$,  then 
$\dh\mbox{\rm proj}_V E<s $ for some $V \in G(n,m)$. In particular, if $E$ satisfies $OSC$ and $0<\dh E <m$ then $\dh \proj_VE <\dh E$ for some $V$.
\end{theo}

A rather different situation occurs if the IFS has {\it dense rotations}, that is the rotation group $G$ is dense in the full group of rotations $SO(n,\mathbb{R})$ or in the group of isometries 
$O(n,\mathbb{R})$. (Note that an IFS of similarities in the plane has dense rotations if  at least one of $O_i$ in \eqref{simform}  is an irrational multiple of $\pi$.) In this case there are no exceptional projections.
\begin{theo}{\rm \cite{FJ,HS,PS}}\label{denserots}
If $E\subset \mathbb{R}^n$ is self-similar with dense rotation group $G$   then 
\begin{equation}\label{hocshm}
\dh\mbox{\rm proj}_V E = \min\{\dh E, m\}\, \mbox{ for  all } V \in G(n,m).
\end{equation}
More generally, $\dh g(E) = \min\{\dh E, m\}$ for all $C^1$ mappings $g: E \to \mathbb{R}^m$ without singular  points, that is the derivative matrix of $g$ everywhere has full rank.
\end{theo}
Note that it is quite  typical for results such as Theorem \ref{denserots} where there are no exceptional directions projections of full dimension to be extendable  by approximation to images under $C^1$ mappings without singular points.

Peres and Shmerkin \cite{PS} proved \eqref{hocshm} in the plane without requiring any separation condition on the IFS. To show this they set up a discrete version of Marstrand's projection theorem to construct a tree of intervals in the subspace (line) $V$  followed by an application of Weyl's equidistribution theorem.  Hochman and Shmerkin \cite{HS}  proved the theorem in higher dimensions, including the extension to $C^1$ mappings, for $E$  satisfying the open set condition. Their proof uses the CP-chains of Furstenberg \cite{Fur,Fur1}, see also \cite{Hoc},  and has three main ingredients: the lower semicontinuity of the expected Hausdorff dimension of the projection of a measure with respect to its `micromeasures', Marstrand's projection theorem, and the invariance of the dimension of projections under the action of the rotation group. Falconer and Jin  \cite{FJ} gave a proof using a compact group extension argument.

 In the dense rotation case, if $\dh E>m$ then  $\dh \mbox{\rm proj}_V E = m$ for all $V \in G(n,m)$ by Theorem \ref{denserots}, but we might also hope that all projections also have positive Lebesgue measure.  Shmerkin and Solomyak \cite{SS} showed that in the plane this is the case for all but a set of directions of Hausdorff dimension 0, see Theorem \ref{projSS}(ii) below, but Rapaport \cite{Rap1} constructed an example with dense rotations in  $\mathbb{R}^2$ such that ${\mathcal L}^1 (\mbox{\rm proj}_\theta E) = 0$ for a dense $G_\delta$ set of $\theta\in [0,\pi)$.
  
The following theorem summarises what is known in the plane, including (iv) and (v) which give packing dimension bounds for the exceptional direcetions, see \cite{Shm} for a detailed discussion.
 
  \begin{theo}\label{projSS}
Let $E$ be the self-similar attractor of an IFS of similarities  \eqref{simform} on $\mathbb{R}^2$.

(i) If $E$ has dense rotations then $\dh \proj_\theta E = \min\{1, \dh E\}$ for all $\theta$.  {\rm \cite{FJ,HS,PS}}.

(ii)  If $\dh E >1$ then $\dh \{\theta: {\mathcal L}^1 (\mbox{\rm proj}_\theta E) = 0\} =0$ {\rm \cite{SS}}.

(iii) If the $a_i$ and the entries in the matrices $r_iO_i$ in the IFS \eqref{simform} are all algebraic numbers then 
$\{\theta: \dh \proj_\theta E < \min\{1, \dh E\} \}$ is countable {\rm \cite{Hoc1}}.

(iv) $\dpk \{\theta: \dh \proj_\theta E < \min\{1, \dh E\} \}=0$ {\rm \cite{Hoc1,Orp}}.

(v) If $0\leq s<\dh E$ then $\dpk \{\theta: \dh \proj_\theta E \leq  s\}\leq s$ {\rm \cite{Orp,Ram}}.
\end{theo}

Recently Algom and Shmerkin \cite{AS} gave a nice criterion that guarantees that the projection of a self-similar $E$ onto a given subspace $V\in G(n,m)$ satisfies Marstrand's theorem; this depends on the curve $\{gV: g\in G\}$ where $G$ is the closure of the group generated by the orthonormal components $O_i$ in the IFS. 

 \subsection{Self-affine sets}
 
If the  IFS consists of affine maps on $\mathbb{R}^n$ of the form 
\begin{equation}\label{affifs}
 f_i(x) = A_i(x) + {t}_i  \quad (1\leq i\leq k),
 \end{equation} 
 where the $ A_i $  are non-singular linear maps with supremum norm $\|A_i\| <1$ and $ t_i$ are translations on $\mathbb{R}^n$, the attractor $E$ defined by \eqref{attractor} is termed {\it self-affine}. (We assume that at least one of the $ A_i $ is not a similarity, otherwise $E$ would be self-similar.)

There is a natural formula for the dimension of $E$ that reflects the distortion of a ball under compositions of the $\{A_i\}$. Recall that the {\it singular values} $\alpha_1\geq\alpha_2 \geq \cdots \geq \alpha_n$ of an $n\times n$ matrix $A$ (or corresponding linear map) are the lengths of the principal semi-axes of the ellipsoid $A(B(0,1))$, equivalently the positive square roots of the eigenvalues of $AA^T$. We define the {\it singular value function} of $A$ for $t>0$ by 
 $$\phi^t(A) = 
\left\{
\begin{array}{lll}
\alpha_1 \cdots  \alpha_k^{t-(k-1)}  &\text{ if }& 0\leq k-1 <t\leq k \leq n      \\
( \alpha_1 \cdots\alpha_d)^{t/n} &  \text{ if }& t\geq n   \\
\end{array}
\right. 
$$
where $k=\lceil t\rceil$ is the least integer greater than or equal to $t$. Then the   {\it affinity dimension} of the IFS attractor $E$ is defined by a pressure-type expression: 
\begin{equation}\label{affdim}
\dim_{\rm aff} E=\dim_{\rm aff} (A_1, \ldots,A_k):=\inf \Big\{t>0:  \sum_{r=1}^\infty  \sum_{1\leq i_1,\ldots,i_r\leq k} \phi^t(A_{i_1}A_{i_2}\cdots A_{i_r})<\infty \Big\}.
\end{equation}
(Strictly, $\dim_{\rm aff} E$ depends on the linear parts of the IFS mappings $A_i$  that define $E$, though this rarely causes confusion.)  
Noting  that for $0<t\leq n$ with $k=\lceil t\rceil$  the singular value function  $\phi^t(A)$ is, to within a bounded multiple, the number of pieces of diameter $\alpha_k$ into which $A(B(0,1))$ can be cut multiplied by $\alpha_k^t$, and using coverings by such pieces, it is easily seen that 
\begin{equation}\label{affdimbound}
\dh E \leq \dlb  E  \leq \dub E \leq \min\{n,  \dim_{\rm aff} E \}. 
\end{equation}
As well as providing this universal upper bound, it has been known for sometime that the affinity dimension equals the Hausdorff and box dimensions for `typical' self-affine sets, in the sense that if $\|A_i\|<\frac12$ for all $i$ (or more generally if $\|A_i\| + \|A_i\|<1$ for all $i\neq j$, see \cite{BSS})  there is equality in \eqref{affdimbound} for $\mathcal{L}^{kn}$-almost all translation vectors $(t_1, \ldots, t_k) \in \mathbb{R}^{kn}$, with the notation of \eqref{attractor} and \eqref{affifs}, see \cite{BSS,Fal7,Fal8,Sol}. Feng, Lo and Ma \cite{FLM} have investigated projections in this context, and obtained natural results for the Hausdorff, packing and box dimensions of $\proj_V E$ for  $\mathcal{L}^{kn}$-almost all translations, with bounds on the Hausdorff dimension of the set of vectors  $(t_1, \ldots, t_k)$ for which the projections have exceptionally small dimensions.

The dimension theory of self-affine sets has recently developed rapidly using techniques from ergodic theory, and equality throughout \eqref{affdimbound} has been established for large classes of  `specific' self-affine sets, starting with $\mathbb{R}^2$, see, for example, \cite{BHR, BSS, MS}.

We call an IFS of affine maps \eqref{affifs} on $\mathbb{R}^2$ {\it strongly irreducible} if there is no finite collection of lines $\{L_j\}$ through the origin  whose union is preserved by all the matrices  $A_i$, that is for all $\{L_j\}$,  $A_i (\bigcup_j L_j) \neq \bigcup_j L_j$ for some $i= 1,\ldots, k$.

\begin{theo}\label{selfaff}
If the self-affine IFS $\{f_i\}_{i=1}^m$ on $\mathbb{R}^2$ satisfies the strong open set condition and is strongly irreducible, then for the attractor $E$,

(i)  $\dh E = \db E = \dim_{\rm aff} E$,

(ii) $\dh \proj_\theta E = \min\{\dh E, 1\}$   for {\emph all } $\theta\in [0,\pi).$
\end{theo}

Theorem \ref{selfaff} was proved by B\'{a}r\'{a}ny, Hochman and Rapaport in 2021 \cite{BHR}, though the projection part (ii) had already been established by Falconer and Kempton in 2017 \cite{FKem} in the case where the matrices $L_i$ have all entries strictly positive. Hochman and Rapaport \cite{HR} have since showed that the strong open set condition can be replaced by the weaker condition of exponential separation.

The 3-dimensional analogue of Theorem \ref{selfaff} has now been established. Rapaport's work \cite{Rap} on self-affine measures  and their projections together with results by Morris and Sert \cite{MSe} on positive subsystems enable Theorem \ref{selfaff} to be extended to IFSs of affine maps on $\mathbb{R}^3$ satisfying the strong separation condition. 

Morris \cite{Mor} has since shown that such projection results are not necessarily valid for self-affine sets in 4-dimensions. He constructs affine IFSs in $\mathbb{R}^4$ that are totally irreducible and satisfy the strong separation condition with attractor $E$ for which there is a 1-parameter family of 2-planes $V_t$ such that $\dh \proj_{V_t} E< \dh E$ for all $t$. The example depends on the fact that the IFS mappings do not act strongly irreducibly  on every exterior power of $\mathbb{R}^4$.

Some remarkable progress has been made very recently on projections of self-affine sets in $\mathbb{R}^n$ for $n \geq 4$. For an affine IFS \eqref{affifs} with $\|A_i\| + \|A_i\|<1$ for all $i\neq j$, in papers posted on the arXiv on the same day, Feng and Xie \cite{FX} and Morris and Sert \cite{MSe2} define a variant $\dim_{\rm aff} (A_1, \ldots,A_k; V)$ of the affinity dimension \eqref{affdim} which additionally depends on the subspace $V\in G(n,m)$ as well as on the $A_i$.  They show that  $\dh \proj_V E= \dh \proj_V E= \dim_{\rm aff} (A_1, \ldots,A_k; V)$ for $\mathcal{L}^{kn}$-almost all translations $(t_1, \ldots, t_k)$. Moreover, $\dim_{\rm aff} (A_1, \ldots,A_k; V)$  can take just finitely many different values for $V\in G(n,m)$; indeed Morris and Sert show that there is a finite filtration  $\varnothing =W_{p+1} \subset W_{p} \subset \cdots \subset W_{1} = G(n,m)$ of algebraic varieties, each invariant under the linear group generated by the $\{A_i\}$, such that $\dh \proj_V E$ is constant for all $V\in W_j \setminus W_{j+1}$.

In a different direction, a recent paper by Lai and Patil \cite{LP} characterises when every projection of a self-affine set coincides with the corresponding projection of its convex hull. Let $\{f_i\}_{i=1}^k$  be an IFS of affine maps on $\mathbb{R}^n$ with attractor $E$ and consider the projections of $E$ onto lines. Suppose that 
$$\bigcup_{i=1}^k f_i({\rm conv}E )= K_1 \cup\cdots \cup K_r$$
where the $\{K_j\}_1^r$ are the disjoint connected components of this set  and ${\rm conv}$ denotes  a convex hull.
 Then $\proj_\theta ({\rm conv}E) = \proj_\theta E$ for all directions $\theta \in S^{n-1}$ if and only if 
for all proper partitious $I\cup I^c$ 
of  $\{1,2,\ldots,r\}$,
$$ {\rm conv}\Big(\bigcup_{i\in I}K_i\Big)  \cap {\rm conv}\Big(\bigcup_{i\in I^c}K_i\Big) \neq \varnothing.$$
This elegant paper contains many examples and ideas related to this notion of `very thick shadows in every direction' for sets that may themselves be very thin.

Iterated function systems with contractions of the form
\begin{equation}
f_i(x,y) = (a_i x +c_i, b_i y +d_i) \label{carpetIFS}
\end{equation}
are called {\it diagonal} IFSs, since the linear parts are represented by diagonal matrices, and the attractors are termed {\it self-affine  carpets}. Such affine transformations leave the horizontal and vertical directions invariant, so these IFSs do not fit into the `irreducible' setting mentioned above.  Well-known examples include Bedford-McMullen, Gatzouras-Lalley or Bara\'{n}ski carpets, see \cite{Fal8,FFJ,Fra4} for definitions and details. Of course, by Marstrand's theorem, $\dh \proj_\theta E = \min\{\dh E,1\}$ for almost all $\theta$, but  
for many self-affine carpets this holds for all directions other than those parallel to the axes, though often there is a requirement that the ratio of the logs of certain of the scaling parameters should be irrational. 
Such results, including for the above carpets, were given by Ferguson, Fraser and Sahlsten  \cite{FJS}. Recently, following work on measures by  Py\"{o}r\"{a}l\"{a} \cite{Pyo}, Feng \cite{FenZ} has shown that the Marstrand conclusion holds for all but the axes directions in a very general setting, namely for diagonal IFSs for which the one-dimensional IFSs induced on each coordinate axis satisfies the exponential separation condition, subject to a weak irrationality condition on the contraction ratios.

Rather little is known about the box and packing dimensions of projections of self-affine carpets. One strange result  has come to light from projection properties of intermediate dimensions. For carpets for which the intermediate dimensions are continuous at 0 (see Section \ref{inter}), which include Bedford-McMullen carpets, if $E \subset \mathbb{R}^2$ with $\dh E <1 \leq \dub E$ then 
$\dub \proj_\theta E < 1$
for all $\theta$. Thus the box-dimensions of the projections of $E$ are restricted by a condition involving the Hausdorff dimension of $E$, see \cite{BFF}.

\subsection{Self-conformal sets}
Let $\{f_i\}_{i=1}^k$ be an IFS of  conformal $C^{1+\epsilon}$-mappings $f_i: U\to U$ on some convex open $U\subset \mathbb{R}^n$ satisfying the open set condition. Thus  $f'_i(x) = r_i(x) O_i(x)$ where $O_i(x)$ are rotations and we assume that $0< r_- \leq r_i(x) \leq r_+ <1$  for $x\in U$; thus the derivatives $f_i^\prime(x)$ are scalar multiples of rotation matrices. The set determined by \eqref{attractor} is then termed a {\it self-conformal} attractor. 

We need a condition that corresponds to `dense rotations' in this non-linear context. We code the points of the attractor in the usual way as $x({\bf i}) = \lim_{p \to \infty} f_{i_1}\circ \cdots \circ f_{i_p}(0)$, where  ${\bf i} = i_1,i_2,\ldots\, (1\leq i_j\leq k)$. We write $G:= SO(n, \mathbb{R})$ for the rotation group
and define $\phi:\{1,2,\ldots,k\}^\mathbb{N}\to G$ by $\phi( {\bf i}) = O_{i_1}(x(\sigma {\bf i}))$, where $\sigma$ is the shift map; thus $\phi( {\bf i})$ is the local rotation at $x(\sigma {\bf i})$ effected by $f_{i_1}^{-1}$. We then define the skew product $\sigma_\phi: \{1,2,\ldots,k\}\times G\ \toitself$ by
\begin{equation}\label{skew}
\sigma_\phi({\bf i},O) = (\sigma {\bf i}, O\phi({\bf i})). 
\end{equation} 
The dynamical system $\sigma_\phi: \{1,2,\ldots,k\}^\mathbb{N}\times G$ is ergodic if  $\sigma_\phi$ has a dense orbit. Using this ergodicity Bruce and Jin \cite{BJ} developed the CP chain and  compact group extension theorem approach in \cite{FJ} to obtain the following  nice theorem.

\begin{theo}\label{conf}
Let $E$ be the attractor of a conformal IFS as above. If $\sigma_\phi$ in \eqref{skew} has a dense orbit then
$$\dh  \proj_V E =  \min\{m,\dh E\}\text{ for all } V\in G(n,m).$$
Moreover, $\dh g(E) = \min\{\dh E, m\}$ for all $C^1$ mappings $g: E \to \mathbb{R}^m$ without singular  points.
\end{theo}
Julia sets of complex mappings provide an important class of self-conformal sets. For  $h: \mathbb{C} \to \mathbb{C}$ given by  $h(z) = z^2 +c$ with  $|c| > \frac14 (5+2\sqrt{6})= 2.475\ldots$, let $f_1,f_2$   be the IFS defined by the two (contracting) branches of $h^{-1}$ on a suitable domain. The self-conformal attractor $E$ of this IFS is precisely the (repelling) Julia set of $h$, see Figure 5. Provided $c$ is such that $\arg (1+ \sqrt{1-4c})/\pi$ is irrational then $\sigma_\phi$ in \eqref{skew} has a dense orbit, so Theorem \ref{conf} implies that all projections have dimension equal to $\dh E$. 

 \begin{figure}[h]
\begin{center}
\includegraphics[width=0.55\linewidth]{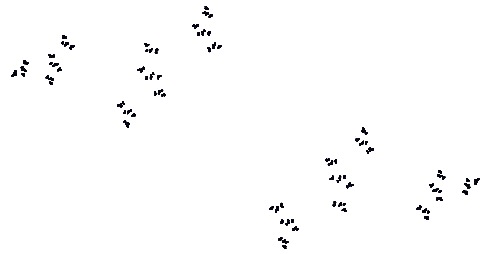}

Figure 5: A self-conformal Julia set.
\end{center}
\end{figure}

\section{Radial projections}
\setcounter{equation}{0}
\setcounter{theo}{0}
\setcounter{figure}{0}
 
Rather than projecting a set  $E \subset \mathbb{R}^n$ onto lines or subspaces one can consider the projection  at  points. The projection at $x \in \mathbb{R}^d$ of $E \subset \mathbb{R}^n$ is the set of directions of the half-lines emanating from $x$ that intersect $E$ (looking at the sky at night one observes a radial projection of the stars). Thus writing $S^{n-1}$ for the  $(n-1)$-dimensional unit sphere embedded in  $\mathbb{R}^n$, we define {\it radial projection at }$x$ to be the mapping   $\proj_x: \mathbb{R}^n\setminus \{x\} \to S^{n-1}$ given by 
\begin{equation*}\label{radproj}
\proj_x(y)  = \frac{y-x}{|y-x|}, \qquad y\in \mathbb{R}^n\setminus \{x\}.
\end{equation*}
[We remark that some properties equivalent to radial projections may be found in Marstrand's 1954 paper \cite{Mar}.]

A Marstrand-type result for radial projections  in the plane follows from Theorem \ref{marthm} by extending the plane to include the line at infinity $L_\infty$ corresponding to directions of parallel lines in the plane. There is a projective transformation $\psi_L$ that maps lines to lines  with any given line $L$ mapped to $L_\infty$ and lines through each point on $L$ mapped to parallel lines in some direction $\theta$. Provided a set $E$ is such that  $\dh (E\setminus L) =\dh E$ it follows that $\dh\psi_L(E) = \dh E$, so Marstrand's theorem gives that $\dh \proj_x E =\min\{\dh E, 1\}$ at $\mathcal{L}$-almost all $x\in L$.  By choosing a suitable set of parallel lines for $L$, it follows using Fubini's theorem that $\dh \proj_x E =\min\{\dh E, 1\}$ for $\mathcal{L}^2$-almost all $x\in \mathbb{R}^2$.

 Bounds on the dimensions of exceptional points $x$ at which the radial projection of $E$ is smaller than the `expected' $\min\{\dh E, n-1\}$ were recently obtained by Orponen, Shmerkin and Wang \cite{OSW} for  the planar case and Bright and Gan  \cite{BG} in higher dimensions.   For  a Borel set  $E\subset \mathbb{R}^n$ with $\dh E \in (k-1,k]$ where $k\in \{1,\ldots, n-1\}$, they obtained  a sharp Kaufman-type bound
\begin{equation}\label{radkau}
 \dh \{x\in  \mathbb{R}^n\setminus E: \dh  \proj_x E < \dh E\} \leq k.
 \end{equation}
For   $\dh E \in (k,k+1]$ where $k\in \{1,\ldots, n-1\}$ and $0\leq s\leq k$,  there is a Falconer-type bound
\begin{equation}\label{radfal}
\dh \{x\in  \mathbb{R}^n\setminus E: \dh   \proj_x E <s \}\leq \max\{k+s-\dh E,0\}.
\end{equation}
The   inequality \eqref{radkau} had been conjectured by \cite{Liu} and \eqref{radfal}  by Lund, Pham and Thu \cite{LPT}. The paper \cite{OSW} provides a good overview and references for earlier work.

{\it Pinned distance maps} are a sort of dual to radial projections. For $x\in \mathbb{R}^n$ we define ${\rm dist}_x: \mathbb{R}^n \to [0,\infty)$ by ${\rm dist}_x(y) = |y-x|$. Thus  for each $x\in \mathbb{R}^n$, ${\rm dist}_x E$ is the aggregate of the distances of points of  $E$ from $x$. The relationships between the dimensions of $E$ and of ${\rm dist}_x E$  are a matter of intense study in relation to the distance set problem. In  one version this asks whether $\dh E>n/2$ implies that $\mathcal{L}(\bigcup_{x\in E}{\rm dist}_x E)\equiv  \mathcal{L}(\{|x-y|:x,y \in E\})>0$. The pinned distance problem is even stronger:  if $\dh E>n/2$ can one find a point $x\in E$ (or indeed many such points) such that $\mathcal{L}({\rm dist}_x E)>0$?  We do not pursue these questions here; there is an enormous literature on these problems, see, for example, \cite{Shm2}.

\section{Some other aspects of fractal projections}\label{lastsec}
\setcounter{equation}{0}
\setcounter{theo}{0}
\setcounter{figure}{0}

There are many other aspects of projections of sets that are actively being researched.   We end with a very brief mention of some of these where recent progress has been made, with references to where more details and further references may be found.
\medskip

\noindent {\it Projections of random sets.} 

Almost any fractal construction may be randomised and if the random process determining the form of the fractals is present at arbitraily small scales the `zero-one law' from probability theory may guarantee that quantities such as the dimensions of the random set and of its projections take a certain value almost surely.

The best-known random fractals are random variants of iterated function system constructions, see Section \ref{secss}. {\it Fractal percolation} is a natural randomisation of the iterated process indicated in \eqref{attint}. 
Here, for a given probability $0<p<1$, the sets $A_{i_1,i_2, \ldots, i_j}$ are independently included in the construction with probability $p$ and omitted with probability $1-p$, to get a (possibly empty) random set $E$ as the intersection of the $j$th level random approximations.

For the most frequently considered case, let $A$ be the unit square divided into $M\times M$ subsquares of side $1/M$ in the natural way, where $M\geq 2$ is an integer. Let the $f_i \, (1\leq i \leq M^2)$ be the homotheties (similarities without rotation or reflection) that map $A$ onto each of these subsquares. This is {\it Mandelbrot percolation} which yields a set $E$ which almost surely has box and Hausdorff dimensions $2+\log p/\log M$, subject to non-extinction. Then, for projections of the Mandelbrot percolation set, almost surely  the conclusions of Marstrand's Theorem \ref{marthm} hold for all projections, with further conditions ensuring that the projections all contain an interval, see  \cite{BF, FFJ,PR,RS,RS2, RS3}.  There are many variants, for example the probability of retaining subsets need not be constant throughout, see for example \cite{RS2}, or for 3-D variants see \cite{SV,OSi}, including percolation on the Menger sponge.

For percolation based on general IFSs of similarities, if the underlying IFS has dense rotations, see Section \ref{sss}, ergodic theoretic methods yield a natural random analogue of Theorem \ref{denserots}  \cite{FJ,FJ2}.  Marstrand-type results for projections of random covering `lim-sup' sets are considered in \cite{CKLS}. Shmerkin and Suomala \cite{SSou} introduced a very general theory showing that for a class of random measures, termed spatially independent martingales, strong results are valid for dimensions of projections of the measures and thus of underlying support sets, with many conclusions almost surely holding  for projections in all directions or onto all subspaces. 
\medskip

\noindent {\it Generalised projections.}

Families of projections parameterised by $\theta \in [0,\pi)$ or $V\in G(n,m)$ are particular cases of more general parameterised families of functions, which may be nonlinear, with  
Marstrand and Mattila's projection theorems particular instances. The potential theoretic and Fourier proofs developed by Kaufman and Falconer (see Sections 2 and 3) depend on what is now known as a transversality condition which can often be applied to nonlinear mappings. This idea was developed by Peres and Schlag \cite{PS} who obtained Marstrand-type theorems for the dimension of the images of a set under parameterised families of perhaps nonlinear mappings valid for almost all parameters, together with bounds on the dimension of the exceptional parameters. For recent treatments of transversality see \cite{BSS,BT,Mat2}. Shmerkin \cite{Shm3} used a discretised approach to extend Bourgain's projection theorem on dimensions of exceptional sets to images of all Borel sets under parameterised families of mappings without singular points.

For images of certain classes of set it may be possible to obtain sure results. 
For example, in Theorem \ref{denserots} we noted that, for self-similar sets $E$ with dense rotations, not only do projections in all directions have Hausdorff dimension $\min\{\dh E,m\}$, but this is also true for images of $E$  under $C^1$ maps without singular points. For sets that are self-similar under IFSs of homotopies,  a similar conclusion is true, see \cite{Bar,BLZ}.

Fraser \cite{Fra5}  extended the inequality \eqref{assproj} to show that  
$${\dim}_{\rm A} g(E) \geq \min\{m, {\dim}_{\rm A} E\}$$
 for all sets $E$ and all $C^1$ mappings $g: E \to \mathbb{R}^m$ without singular points, a result with a number of nice corollaries, including a proof of the `distance conjecture' for Assouad dimension.

\medskip

\noindent {\it Projections in Heisenberg groups.}   

The Heisenberg group $\mathbb{H}$ may be identified with $\mathbb{C}\times \mathbb{R} \equiv \mathbb{R}^3$, with the group product given by $(z,t)\ast (w,u) = (z+w, t+u+2\,{\rm Im}(z{\overline w}))$. 
The Heisenberg metric is defined as $d_\mathbb{H}((z,t),(w,u))= \| (w,u)^{-1}\ast(z,t)\|_\mathbb{H}$, where the norm $\| \cdot\|_\mathbb{H}$ is given by $\|(z,t)\|_\mathbb{H}= (|z|^4 +t^2)^{1/2}$.
Then $d_\mathbb{H}$ is left invariant under the group action and behaves very differently  from the Euclidean metric, with the Hausdorff dimension of a subset of $\mathbb{H}$ depending on which metric is used. Despite the lack of isotropy, there is enough geometric structure for projections onto certain families of subspaces to have interesting properties. In particular, for each $\theta\in [0,\pi)$ there is a semidirect group splitting $\mathbb{H}= W_\theta \ast P_\theta$ where $V_\theta$ is a `horizontal' line and $W_\theta$ a 2-dimensional `vertical plane' so we may consider projections from $\mathbb{H}$ to $W_\theta$ and from $\mathbb{H}$ to $P_\theta$ for each $\theta$ corresponding to this splitting. 
Various estimates for the dimension of projections of a Borel set $E\subset \mathbb{H}$ in terms of the dimension of $E$ for almost all $\theta$ have been obtained for both horizontal projections and vertical projections, see \cite{BDFMT, BFMT,FHov,FO2,Har4, Har3,Mat5}, along with analogues for the higher order Heisenberg groups $\mathbb{H}^n$.
\medskip

\noindent {\it Projections in other spaces.}

Balogh and Iseli  have  obtained Marstrand-type results for orthogonal projections defined by geodesics on Riemann surfaces of constant curvature \cite{BI3}. They have also obtained projection theorems in $\mathbb{R}^2$ under general norms, where there is no notion of orthogonality, by replacing projections by nearest point mappings onto lines \cite{BI}. Similarly, in hyperbolic spaces they consider the nearest point map onto images of $m$-planes under an exponential map \cite{BI2}.  Versions for M\"{o}bius transformations on the Riemann sphere and for projective transformations on the projective plane are presented in \cite{IL}.
The proofs in these different settings use the general transversality framework of Peres-Schlag \cite{PS} to obtain projection theorems along with bounds for the exceptional sets of projections.

Linear embeddings of attractors in infinite dimensional dynamical systems into finite dimensional spaces may be regarded as projections. A classical dimension result is due to Hunt and Kaloshin \cite{HK}. Let $E$ be a compact subset of a Banach space $X$ with box-counting dimension $d$. Then for almost every projection or bounded linear function $\pi: X \to \mathbb{R}^m$ such that $m >2d$,
$$\frac{m-2d}{m(1+d)} \dh E \leq \dh \pi(E) \leq \dh E.$$
Here `almost every' is interpreted in the sense of {\it prevalence}, which is a measure-theoretic way of defining sparse and full sets for infinite-dimensional spaces. Further references to projections from this embedding viewpoint include \cite{MR,OHK, Rob,Rob2,RoS}. This is an area where Assouad dimension and the notion of 'thickness' play important roles.

\medskip

\noindent {\it Visible parts of sets.} 

The {\it visible part} $\mbox{Vis}_\theta E$ of a compact set $E\subset \mathbb{R}^2$ from direction $\theta$ is
the set of $x \in E$  such that the half-line from $x$ in direction $\theta$ intersects $E$ in the single point
$x$; thus $\mbox{Vis}_\theta E$ may be thought of as the part of $E$ that can be `seen from infinity' in direction $\theta$. It is immediate from Marstrand's Theorem \ref{marthm} that, for almost all $\theta$, 
$$\dh  \mbox{Vis}_\theta E = \dh E \mbox{ if } \dh E \leq 1 \quad \mbox{ and } \quad \dh \mbox{Vis}_\theta E \geq 1 \mbox{ if } \dh E \geq 1.$$
 It has been long conjectured that if $\dh E \geq 1$  then $\dh  \mbox{Vis}_\theta E = 1$ for almost all $\theta$, but this has only been established for certain specific classes of $E$.The conjecture is easily verified if $E$ is the graph of a function (the only possible exceptional direction being perpendicular to the $x$-axis), see \cite{JJMO}. It is also true for quasi-circles \cite{JJMO} and almost surely for Mandelbrot percolation sets \cite{AJJRS}. For self-similar sets the conjecture is valid if the rotation group is finite and the projection is a countable union of intervals \cite{JJSW}, or if $E$ satisfies the open set condition for a convex open set such that $\proj_\theta E$ is an interval for all $\theta$ \cite{FF} (in this case $E$ need not be connected), as well as for certain self-affine sets \cite{Ros}.  
 
 O'Neil \cite{On} obtained an upper bound for the typical Hausdorff dimensions of the visible sets in terms of $\dh E$.  The analogous conjecture in higher dimensions, that the dimension of the visible part of a compact $E\subset \mathbb{R}^n$ equals $\min\{\dh E, n-1\}$, is also unresolved, but Orponen \cite{Orp8} showed that $\dh  \mbox{Vis}_\theta E \leq n -1/50n$ for almost all $\theta$ for all $n\geq 2$, and very recently D\c{a}browski \cite{Dar} improved this to $\dh   \mbox{Vis}_\theta E  \leq n -1/6$, leading to improved bounds for Ahlfors regular sets.

  \medskip

\section{Final remarks}\label{fin}
\setcounter{equation}{0}
\setcounter{theo}{0}
\setcounter{figure}{0}

This article has surveyed some of the numerous results that may be viewed as descendants of Marstrand's projection theorems. In recent years, fractal geometry has become an established area of mathematics in its own right attracting highly talented mathematicians. Since the prequel survey \cite{FFJ} ten years ago questions have been resolved that even then were considered intractable. 

So many papers have now been written on a large range of aspects of fractal projections that it is impossible to mention them all or include a complete history of the topics highlighted, though papers often review relevant history in their introductions and contain many further references. I apologise to all those whose deserving work has not specifically been mentioned here.

\begin{acknowledgement}
I am very grateful to the many friends and colleagues with whom I have, over the years, discussed fractals and in particular their projections. I thank Jonathan Fraser, Tuomas Orponen and a referee for their comments on drafts of this survey.
\end{acknowledgement}


\begin{thebibliography}{99}

\bibitem{AS} A. Algom and P. Shmerkin. 
On the dimension of orthogonal projections of self-similar measures, 
{\it J. Lond. Math. Soc.(2)} {\bf 112} (2025), Paper No. e70245.

\bibitem{AJJRS}
I. Arhosalo, E. J\"{a}rvenp\"{a}\"{a}, M. J\"{a}rvenp\"{a}\"{a}, M. Rams and P. Shmerkin. Visible parts of fractal percolation, {\it  Proc. Edinb. Math. Soc.} {\bf 55} (2012),  311--331.

\bibitem{BDFMT}
Z.M. Balogh, E. Durand-Cartagena, K. F\"{a}ssler, P. Mattila and J.T. Tyson. The effect of projections on dimension in the Heisenberg group, 
{\it Rev. Mat. Iberoam.} {\bf 29} (2013), 381--432.

\bibitem{BFMT}
Z.M. Balogh, K. F\"{a}ssler, P. Mattila and J.T. Tyson. Projection and slicing theorems in Heisenberg groups, 
{\it Adv. Math.} {\bf 231} (2012), 569--604.

\bibitem{BI3}
Z.M. Balogh and A. Iseli. Dimensions of projections of sets on Riemannian surfaces of constant curvature, 
{\it Proc. Amer. Math. Soc.} {\bf 144} (2016), 2939--2951.

\bibitem{BI}
Z.M. Balogh and A. Iseli. Marstrand type projection theorems for normed spaces, 
{\it J. Fractal Geom.} {\bf 4} (2019), 367--392.

\bibitem{BI2}
Z.M. Balogh and A. Iseli. Projection theorems in hyperbolic space, 
{\it Arch. Math. (Basel)} {\bf 112} (2019), 329--336.

\bibitem{BFZ} C. Bandt, K. Falconer and M. Z\"{a}hle (eds.). {\em Fractal Geometry and Stochastics V},  
{\it  Progress in Probability} {\bf 70}, Birkhauser, Cham, 2015.

\bibitem{Bar} B. B\'{a}r\'{a}ny. 
On some non-linear projections of self-similar sets in $\mathbb{R}^3$, {\it Fund. Math.} {\bf 237} (2017), 83--100. 

\bibitem{BHR} B. B\'{a}r\'{a}ny, M. Hochman, and A. Rapaport. 
Hausdorff dimension of planar self-affine sets and measures, {\it Invent. Math.} {\bf 216} (2019), 601--659. 

\bibitem{BSS} B. B\'{a}r\'{a}ny, K. Simon and B. Solomyak. 
{\it Self-Similar and Self-Affine Sets and Measures}, Princeton U.P., 2024.

\bibitem{BF}
J. Barral and D.-J. Feng. Projection of planar Mandelbrot measures, 
{\it Adv. Math.} {\bf 325} (2018), 640--718.

\bibitem{BLZ}
M. Bond, I. \L{}aba and J. Zahl.  Quantitative visibility estimates for unrectifiable sets in the plane, {\it  Trans. Amer. Math. Soc.} {\bf 368} (2016), 5475--5513.

\bibitem{BT}
R. Bongers and K. Taylor.  Transversal families of nonlinear projections and generalizations of Favard length, {\it  Anal. PDE.} {\bf 16} (2023), 279--308.

\bibitem{Bou}
J. Bourgain.  On the Erd\H{o}s-Volkmann and Katz-Tao ring conjectures, 
{\it  Geom. Funct. Anal.} {\bf 13} (2003), 334--365.

\bibitem{Bou2}
J. Bourgain.  The discretized sum-product and projection theorems, 
{\it  J. Anal. Math.} {\bf 112} (2010), 193--236.

\bibitem{BG}
C. Bright  and S. Gan.
Exceptional set estimates for radial projections in $\mathbb{R}^n$,
 {\em  Ann. Acad. Sci. Fenn. Math.} {\bf 49} (2024), 631--661.
 
\bibitem{BJ}
C. Bruce  and X. Jin.
Projections of Gibbs measures on self-conformal sets,
 {\em  Nonlinearity} {\bf 32} (2019), 603--621.
 
\bibitem{BFF}
S. Burrell, K.J. Falconer and J. Fraser.
Projection theorems for intermediate dimensions,
 {\em  J. Fractal Geom.} {\bf 8} (2021), 95--116.

\bibitem{CKLS}
C. Chen, H. Koivusalo, B. Li, Bing and V. Suomala. 
Projections of random covering sets,
 {\em  J. Fractal Geom.} {\bf 1} (2014), 449--467.

\bibitem{Dar}
D. D\c{a}browski.
Visible parts and slices of Ahlfors regular sets,
 {\em  Discrete Anal. 2024},  Paper No. 17.

\bibitem{EM}
G. Edgar and C. Miller. Borel subrings of the reals, 
 {\em  Proc. Amer. Math. Soc.} {\bf 131} (2003), 1121--1129.
 
\bibitem{Fal1}
K.J. Falconer. Hausdorff dimension and the exceptional set of projections, 
{\it Mathematika} {\bf 29} (1982), 109--115.

\bibitem{Fal7}
K.J. Falconer. The Hausdorff dimension of self-affine fractals, 
{\it Math. Proc. Cambridge Philos. Soc.} {\bf 103} (1988), 339--350.

\bibitem{Fal8}
K.J. Falconer, Dimensions of Self-affine Sets: A Survey,
in {\em Further Developments in Fractals and Related Fields},  J. Barral and S. Seuret (eds.), pp 115--134,
{\it  Trends in Mathematics}, Springer, New York, 2013.

\bibitem{Fal}
K.J. Falconer.
{\em Fractal Geometry: Mathematical Foundations and Applications},
 John Wiley \& Sons, Hoboken, NJ, 3rd. ed., 2014.
 
\bibitem{Fal5}
K.J. Falconer. A capacity approach to box and packing dimensions of projections and other images, 
in {\em Analysis, Probability and Mathematical Physics on Fractals},  pp 1--19,
{\it  Fractals Dyn. Math. Sci. Arts Theory Appl.} {\bf 5}, World Sci. Publ., Hackensack, 2020.
 
\bibitem{Fal4}
K.J. Falconer. A capacity approach to box and packing dimensions of projections of sets and exceptional directions, {\it J. Fractal Geom.} {\bf 8} (2021), 1--26. 

\bibitem{Fal6}
K.J. Falconer. Intermediate Dimensions: A Survey, 
in {\em Thermodynamic Formalism}, pp 469--493, M. Pollicott and S. Vaienti (eds)
{\it  Lecture Notes in Mathematics} {\bf 2290}, Springer, Cham, 2021.

\bibitem{FF}
K.J. Falconer and J.M. Fraser. The visible part of plane self-similar sets, 
{\it Proc. Amer. Math. Soc.} {\bf 141} (2013), 269--278.
\bibitem{FFJ}
K.J. Falconer, J.M. Fraser and X. Jin. Sixty years of fractal projections, 
in {\em Fractal Geometry and Stochastics V}, C. Bandt, K. Falconer and M. Z\"{a}hle (eds.), pp 3--25,
{\it  Progress in Probability} {\bf 70}, Birkhauser, Cham, 2015.

\bibitem{FFK}
K.J. Falconer, J.M. Fraser and T. Kempton. Intermediate dimensions. {\it Math. Z.} 
{\bf 296} (2020), 813--830.

\bibitem{FH}
K.J. Falconer and J.D. Howroyd. Projection theorems for box and packing dimensions, 
{\it Math. Proc. Cambridge Philos. Soc.} {\bf 119} (1996), 287--295.

\bibitem{FH2}
K.J. Falconer and J.D. Howroyd. Packing dimensions of projections
and dimension profiles, {\it Math. Proc. Cambridge Philos. Soc.} {\bf 121} (1997), 269--286.

\bibitem{FJ}
K.~J. Falconer and X. Jin.
Exact dimensionality and projections of random self-similar measures and sets,
{\em J. Lond. Math. Soc.(2)} {\bf 90} (2014), 388--412.

\bibitem{FJ2}
K.~J. Falconer and X. Jin.
Self-similar sets: projections, sections and percolation, in
{\em Recent Developments in Fractals and Related Fields},  J. Barral and S. Seuret (eds.), pp. 113--127, Springer, Cham, 2015.

\bibitem{FKem}
K.~J. Falconer and T. Kempton.
The dimension of projections of self-affine sets and measures,
{\it Ann. Acad. Sci. Fenn. Math.} {\bf 42} (2017), 473--486.

\bibitem{FM}
K.J. Falconer and P. Mattila. The packing dimension of projections and sections of
measures, {\it Math. Proc. Cambridge Philos. Soc.} {\bf 119} (1996), 695--713.

\bibitem{Far}
A. Farkas.
Projections of self-similar sets with no separation condition, {\it Israel J. Math.} {\bf 214} (2016), 67--107.

\bibitem{FHov}
K. F\"{a}ssler and R. Hovila. Improved Hausdorff dimension estimates for vertical projections in the Heisenberg group, {\it Ann. Sc. Norm. Super. Pisa. Cl. Sci (5)} {\bf 15} (2016), 459--483.

\bibitem{FO1}
K. F\"{a}ssler and T. Orponen. On restricted families of projections in $\mathbb{R}^3$, {\it Proc. London Math. Soc.(3)} {\bf 109} (2014), 353--381.

\bibitem{FO2}
K. F\"{a}ssler and T. Orponen. Vertical projections in the Heisenberg group via cinematic functions and point-plate incidences, {\it Adv. Math.} {\bf 431} (2023), Paper No. 109248.

\bibitem{FLM} D.-J. Feng, C.-H Lo and C.-Y. Ma.
Dimensions of projected sets and measures on typical self-affine sets,
{\em Adv. Math} {\bf 431} (2023), Paper No. 109237.

\bibitem{FX} D.-J. Feng and Y.-H. Xie. 
Dimensions of orthogonal projections of typical self-affine sets and measures,  arXiv:2502.04000 (2025).

\bibitem{FenZ} Z. Feng. 
Dimension of diagonal self-affine measures with exponentially separated projections,  arXiv:2501.17378 v.2 (2025).

\bibitem{FJS}
A. Ferguson, T. Jordan, P. Shmerkin.
The Hausdorff dimension of the projections of self-affine carpets,
{\em Fund. Math.} {\bf 209} (2010), 193--213.

 \bibitem{Fra1}
J.M. Fraser. 
Distance sets, orthogonal projections, and passing to weak tangents,
{\em Israel J. Math.} {\bf 226} (2018), 851--875.

 \bibitem{FraBk}
J.M. Fraser. {\it Assouad Dimension and Fractal Geometry},
Cambridge University Press, Cambridge, 2021.

\bibitem{Fra4}
J.M. Fraser. Fractal geometry of Bedford-McMullen carpets, 
in {\em Thermodynamic Formalism}, 495--516,
{\it Lecture Notes in Math.} {\bf 2290}, Springer, Cham, 2021.

 \bibitem{Fra5}
J.M. Fraser. 
A nonlinear projection theorem for Assouad dimension and applications,
{\em  J. Lond. Math. Soc.(3)}  {\bf 107} (2023), 777--797.

 \bibitem{Fra3}
J.M. Fraser. 
The Fourier spectrum and sumset type problems,
{\em  Math. Ann.}  {\bf 390} (2024), 3891--3930.

 \bibitem{Fra2}
J.M. Fraser. 
Applications of dimension interpolation to orthogonal projections,
{\em Res. Math. Sci.}{\bf 12} (2025), no. 1, Paper No. 10.

  \bibitem{FK}
J.M. Fraser and A. K\"{a}enm\"{a}ki. 
Attainable values for the Assouad dimension of projections,
{\em Proc. Amer. Math. Soc.} {\bf 148} (2020), 3393--3405.

\bibitem{FO}
J.M. Fraser and A.E. deOrellana. 
A Fourier analytic approach to exceptional set estimates for orthogonal projections,
{\em Indiana Univ. Math. J.} to appear, arXiv:2404.11179 (2024).

 \bibitem{FdO}
J.M. Fraser and T. Orponen. 
The Assouad dimension of projections of planar sets,
{\em Proc. Lond. Math. Soc.(2)} {\bf 114} (2017), 374--398.


 \bibitem{Fur}
H. Furstenberg.
Ergodic fractal measures and dimension conservation,
{\em Ergodic Theory Dynam. Systems} {\bf 28} (2008), 405--422.

 \bibitem{Fur1}
H. Furstenberg. 
Ergodic theory and fractal geometry, {\it CBMS Regional Conference Series in Mathematics} {\bf 120} (2014),
American Mathematical Society. Providence, RI.

 \bibitem{GGGHMW}
S. Gan, S. Guo, L. Guth, T.L.J. Harris, D. Maldague and H. Wang.
On restricted projections to planes in $\mathbb{R}^3$,
arXiv:2207.13844 (2022). 

 \bibitem{GGW}
S. Gan, S. Guo and  H. Wang.
A restricted projection problem for fractal sets in $\mathbb{R}^n$,
{\it Camb. J. Math.}  {\bf 12} (2024), 535--561.

 \bibitem{GGM}
S. Gan, L. Guth and D. Maldague.
An exceptional set estimate for restricted projections to lines in $\mathbb{R}^3$,
{\em J. Geom. Anal.} {\bf 34} (2024), no. 1, Paper No. 15. 

 \bibitem{Har2}
 T.L.J. Harris. Improved bounds for restricted projection families via weighted Fourier restriction, {\it  Ann. Acad. Sci. Fenn. Math.} {\bf 45} (2020), 723--737.
 
\bibitem{Har4}
 T.L.J. Harris. An a.e. lower bound for Hausdorff dimension under vertical projections in the Heisenberg group, {\it  Anal. PDE} {\bf 15} (2022), 1655--1701.
 
 \bibitem{Har3}
 T.L.J. Harris. A Euclidean Fourier-analytic approach to vertical projections in the Heisenberg group, {\it  Bull. Lond. Math. Soc.} {\bf 55} (2023), 961--977.

\bibitem{Har}
 T.L.J. Harris. Length of sets under restricted families of projections onto lines, in {\em Recent Developments in Harmonic Analysis and its Applications Contemporary Mathematics,} 1-17, {\em Contemp. Math. }{\bf  792}, Amer. Math. Soc., Providence,  2024.
 
 \bibitem{Har1}
 T.L.J. Harris. Exceptional sets for length  under restricted families of projections onto lines in $\mathbb{R}^3$, arXiv:2408.04885 (2025).
 
\bibitem{He}
 W. He. Orthogonal projections of discretized sets, {\em J. Fractal Geom.} {\bf  7}  (2020), 271--317.

\bibitem{Hoc}
M.~Hochman.
Dynamics on fractals and fractal distributions, ar{X}iv:1008.3731v2 (2013).

 \bibitem{Hoc1}
M.~Hochman.
On self-similar sets with overlaps and inverse theorems for entropy,
{\em Ann. of Math.(2)} {\bf 180} (2014), 773--822.

 \bibitem{HR}
 M. Hochman and A. Rapaport. Hausdorff dimension of planar self-affine sets and measures with overlaps, {\em J. Eur. Math. Soc.} {\bf  24}  (2022), 2361--2441.

\bibitem{HS}
M.~Hochman and P.~Shmerkin.
Local entropy averages and projections of fractal measures,
{\em Ann. of Math.(2)} {\bf 175} (2012), 1001--1059.

\bibitem{HK}
B.R. Hunt and Y. Kaloshin. Regularity of embeddings of infinite-dimensional fractal sets into finite-dimensional spaces,
{\it Nonlinearity} {\bf 12} (1999), 1263--1275.

\bibitem{Hut}
J.E. Hutchinson. Fractals and self-similarity, 
{\it Indiana Univ. Math. J.} {\bf 30} (1981), 713--747.

\bibitem{IL}
A. Iseli and A. Lukyanenko. Projection theorems for linear-fractional families of projections, 
{\it Math. Proc. Cambridge Philos. Soc.} {\bf 175} (2023), 625--647.

\bibitem{JJK}
E. J\"{a}rvanp\"{a}\"{a}, M. J\"{a}rvanp\"{a}\"{a} and T. Keleti, Hausdorff dimension and non-degenerate
families of projections, 
{\it J. Geom. Anal.}  {\bf 24} (2014), 2024--2034.

\bibitem{JJLL}
E. J\"{a}rvanp\"{a}\"{a}, M. J\"{a}rvanp\"{a}\"{a}, F. Ledrappier and M. Leikas. One-dimensional families of projections, 
{\it Nonlinearity} {\bf 21} (2008), 453--463.

\bibitem{JJMO}
E. J\"{a}rvanp\"{a}\"{a}, M. J\"{a}rvanp\"{a}\"{a}, P. MacManus and T.C. O'Neil. Visible parts and dimensions, 
{\it Nonlinearity} {\bf 16} (2003), 803--818.

\bibitem{JJSW}
E. J\"{a}rvanp\"{a}\"{a}, M. J\"{a}rvanp\"{a}\"{a}, V. Suomala and M. Wu. On dimensions of visible parts of self-similar sets with finite rotation group, 
{\it Proc. Amer. Math. Soc} {\bf 150} (2022), 2983--2995.

\bibitem{Jar}
M. J\"{a}rvanp\"{a}\"{a}. On the upper Minkowski dimension, the packing dimension, and othogonal
projections, 
{\it Ann. Acad. Sci. Fenn. A Dissertat.} {\bf 99} (1994).

\bibitem{KOV}
A. K\"{a}enm\"{a}ki, T. Orponen and L. Venieri.
A Marstrand-type restricted projection theorem in $\mathbb{R}^3$,
{\it  Amer. J. Math.} {\bf 147} (2025), 81--123.

\bibitem{Kau}
R. Kaufman.
On Hausdorff dimension of projections,
 {\em Mathematika} {\bf 15} (1968), 153-155.

\bibitem{KX}
D. Khoshnevisan and Xiao.
Packing-dimension profiles and fractional Brownian motion,
 {\em Math. Proc. Cambridge Phil. Soc.} {\bf 145} (2008), 145--213.
 
 \bibitem{LP}
C.K. Lai and L.P. Patil.
 Projections of totally disconnected thin fractals with very thick shadows on $\mathbb{R}^d$,
 {\em J. Fractal Geom.} {\bf 11} (2024),  289-317.

\bibitem{LMWY}
E. Lindenstrauss, A. Mohammadi, Z. Wang and L. Yang. An effective version of the Oppenheim conjecture with a polynomial error rate, arXiv:2305.18271 (2023).

\bibitem{Liu}
B. Liu. On Hausdorff dimension of radial projections, 
{\it Rev. Mat. Iberoam.} {\bf 37} (2020), 1307--1319.

\bibitem{LPT}
B. Lund,  T. Pham and V.T.H. Thu. Radial projection theorems in finite spaces,
 {\em Proc. Amer. Math. Soc.}  {\bf 153} (2025), 2169--2183.

\bibitem{LS}
N. Lutz and D. M. Stull.
 Projection theorems using effective dimension,
 {\em Inform. and Comput.} {\bf 297} (2024), Paper No. 105137.

\bibitem{MR}
A. Margaris and J.C. Robinson.
 Embedding properties of sets with finite box-counting dimension,
 {\em Nonlinearity} {\bf 32} (2019), 3523--3547.

\bibitem{Marg}
G.A. Margulis.
 Discrete subgroups and ergodic theory,
 {\it Number theory, trace formulas and discrete groups (Oslo, 1987), 377--398},  Academic Press, Boston, 1989. 
 
\bibitem{Mar}
J.M. Marstrand.
 Some fundamental geometrical properties of plane sets of fractional
  dimensions, {\em Proc. London Math. Soc.(3)} {\bf 4} (1954), 257--302.

\bibitem{Mat}
P.~Mattila.
 Hausdorff dimension, orthogonal projections and intersections with
  planes,
 {\em Ann. Acad. Sci. Fenn. A Math.} {\bf 1} (1975),  227--244.

\bibitem{Mat1}
P. Mattila.
{\em Geometry of Sets and Measures in {E}uclidean Spaces},
Cambridge University Press, Cambridge, 1995.

\bibitem{Mat3}
P. Mattila.
 Hausdorff  dimension, projections, and the Fourier transform,
{\it  Publ. Mat.} {\bf 48} (2004), 3--48.

\bibitem{Mat4}
P. Mattila. Marstrand's Theorems,
in {\em All That Math, Portraits of Mathematicians as Young Readers},  A. C\'{o}rdoba, J.L. F\'{e}rnandez and P.  F\'{e}rnandez (eds.), pp 283--301, {\it  Rev. Mat. Iberoam.}, 2011.

\bibitem{Mat5}
P. Mattila. Recent progress on dimensions of projections,
in {\em Geometry and Analysis of Fractals}, D.-J. Feng and K.-S. Lau (eds.), pp 283--301,
{\it  Springer Proceedings in Mathematics \& Statistics} {\bf 88}, Springer, Heidelberg, 2014.

\bibitem{Mat2}
P. Mattila.
{\em Fourier Analysis and Hausdorff Dimension},
Cambridge University Press, Cambridge, 2015.

\bibitem{Mat6}
P. Mattila. Hausdorff dimension, projections, intersections, and Besicovitch sets,
in {\em New Trends in Applied Harmonic Analysis, Vol. 2 -- harmonic analysis, geometric measure theory, and applications}, pp 129--157,
{\it   Appl. Numer. Harmon. Anal.} {\bf 88}, Birkh\"{a}user/Springer, Cham, 2019.


\bibitem{Mat7}
P. Mattila.
Hausdorff dimension and projections related to intersections,
{\it  Publ. Mat.} {\bf 66} (2022), 305--323.

\bibitem{Mor} I.D. Morris. 
Exceptional projections of self-affine sets with strong irreducibility, arXiv:2308.04894 (2023).

\bibitem{MSe} I.D. Morris and C. Sert. 
A variational principle relating self-affine measures to self-affine sets,  arXiv:2303.03437 (2023).

\bibitem{MSe2} I.D. Morris and C. Sert. 
Projections of self-affine fractals,  arXiv:2502.04001 (2025).

\bibitem{MS} I.D. Morris and P. Shmerkin. 
On equality of Hausdorff and affinity dimensions, via self-affine measures on positive subsystems, 
{\it Trans. Amer. Math. Soc.} {\bf  371} (2019), 1547--1582.

\bibitem{Obe}
D.M. Oberlin.  Restricted Radon transforms and projections of planar sets, 
{\it  Canad. Math. Bull.} {\bf 55} (2012), 815--820.

\bibitem{Obe1}
D.M. Oberlin. Exceptional sets of projections, unions of $k$-planes and associated transforms, 
{\it  Israel J. Math.} {\bf 202} (2014), 331--342.

\bibitem{OO}
D.M. Oberlin and R. Oberlin. Application of a Fourier restriction theorem to certain families
of projections in $\mathbb{R}^3$, 
{\it J. Geom. Anal.} {\bf 25} (2015), 1476--1491.

\bibitem{On}
T.C. O'Neil. The Hausdorff dimension of visible sets of plane continua,  
{\it  Trans. Amer. Math. Soc.} {\bf 359} (2007), 5141--5170.

\bibitem{OSi}
V. Orgov\'{a}nyi and K. Simon. Projections of the random Menger sponge, {\it Asian J. Math.} {\bf 27} (2023), 893--935.

\bibitem{Orp}
T. Orponen. On the packing dimension and category of exceptional sets of orthogonal projections,   {\it Ann. Mat. Pura Appl. (4)} {\bf 194} (2015), 843--880.

\bibitem{Orp1}
T. Orponen. Hausdorff dimension estimates for restricted families of projections in $\mathbb{R}^3$,   {\it Adv. Math.} {\bf 275} (2015),147--183.

\bibitem{Orp7}
T. Orponen. An improved bound on the packing dimension of Furstenberg sets in the plane,  {\it J. Eur. Math. Soc.} {\bf 22} (2020), 797--831.

\bibitem{Orp6}
T. Orponen. On the Assouad dimension of projections,  {\it Proc. Lond. Math. Soc.(3)} {\bf 122} (2021), 317--351.

\bibitem{Orp8}
T. Orponen. On the dimension of visible parts,  {\it J. Eur. Math. Soc.} {\bf 252} (2023), 1969--1983.

\bibitem{Orp5}
T. Orponen. On the discretised ABC sum-product problem,  {\it Trans. Amer. Math. Soc.} {\bf 377} (2024), 4647--4702.

\bibitem{OS}
T. Orponen and P. Shmerkin. On the Hausdorff dimension of Furstenberg sets and orthogonal projections in the plane,  {\it Duke Math. J.}  {\bf 172} (2023), 3559--3632.

\bibitem{OS1}
T. Orponen and P. Shmerkin. Projections, Furstenberg sets, and the {\it abc} sum-product problem, arXiv:2301.10199 (2023).

\bibitem{OSW}
T. Orponen, P. Shmerkin and Hong Wang. Kaufman and Falconer estimates for radial projections and a continuum version of Beck's theorem,  {\it Geom. Funct. Anal.}  {\bf 34} (2024), 164--201.

\bibitem{OV}
T. Orponen and L. Venieri. Improved bounds for restricted families of projections to planes in $\mathbb{R}^3$,  {\it Int. Math. Res. Not. IMRN} {\bf 19} (2020), 5797--5813.

\bibitem{OHK}
W. Ott, B. Hunt and V. Kaloshin. The effect of projections on fractal sets and measures in Banach spaces,  {\it Ergodic Theory Dynam. Systems } {\bf 26} (2006), 869--891.

\bibitem{PR}
Y. Peres and M. Rams. Projections of the natural measure for percolation fractals, {\it Israel J. Math.} {\bf 214} (2016), 539--552.

\bibitem{PSc}
Y. Peres and B. Schlag. Smoothness of projections, Bernoulli convolutions, and the dimension of exceptions, {\it Duke Math. J.} {\bf 102} (2000), 193--251.

\bibitem{PS}
Y.~Peres and P.~Shmerkin.
Resonance between Cantor sets,
 {\em Ergodic Theory Dynam. Systems.} {\bf 29} (2009), 201--221.
 
\bibitem{PYZ}
M. Pramanik, T. Yang and J. Zahl.
 A Furstenberg-type problem for circles, and a Kaufman-type restricted projection theorem in $\mathbb{R}^3$,
 arXiv:2207.02259 v.3 (2024).

\bibitem{Pyo}
A. Py\"{o}r\"{a}l\"{a}.
The dimension of projections of planar diagonal self-affine measures,
 {\em Ann. Acad. Sci. Fenn. Math.}, {\bf 50} (2025), 59--78.
 
\bibitem{Ram}
M.~Rams.
 Packing dimension estimation for exceptional parameters,
 {\em Israel J. Math.} {\bf 130} (2014), 125--144.

\bibitem{RS}
M.~Rams and K.~Simon.
 The dimension of projections of fractal percolations,
 {\em J. Stat. Phys.} {\bf 154} (2014), 633--655.

\bibitem{RS3}
M.~Rams and K.~Simon.
The geometry of fractal percolation,
in {\em Geometry and Analysis of Fractals} D.-J. Feng and K.-S. Lau (eds.), pp 303--324,
{\it  Springer Proceedings in Mathematics \& Statistics.} {\bf 88}, Springer-Verlag, Berlin Heidelberg, 2014. 

\bibitem{RS2}
M.~Rams and K.~Simon.
Projections of fractal percolations,
{\em Ergodic Theory Dynam. Systems} {\bf  35} (2015),  530--545.

\bibitem{Rap1} A. Rapaport. A self-similar measure with dense rotations, singular projections and discrete slices,  
{\it Adv. Math.} {\bf  321} (2017),  529--546. 

\bibitem{Rap} A. Rapaport. On self-affine measures associated to strongly irreducible and proximal systems,  {\it Adv. Math.} {\bf  449} (2024),  Paper No. 109734. 

\bibitem{RW}
K. Ren and H. Wang.
 Furstenberg sets estimates in the plane,
  arXiv:2308.08819  (2023).
  
\bibitem{Rob2} J.C. Robinson. Linear embeddings of finite-dimensional subsets of Banach spaces into Euclidean spaces,  
{\it Nonlinearity} {\bf  22} (2009),  711--728. 

\bibitem{Rob}
J.C. Robinson.
{\em Dimensions, Embeddings, and Attractors},
Cambridge University Press, Cambridge, 2010.

\bibitem{Ros} E. Rossi. Visible part of dominated self-affine sets in the plane,  
{\it Ann. Acad. Sci. Fenn. Math.} {\bf  46} (2021),  1089--1103. 

\bibitem{RoS} E. Rossi and P. Shmerkin. H\"{o}lder coverings of sets of small dimension,  
{\it J. Fractal Geom.} {\bf  6} (2019),  285--299. 

\bibitem{Shm}
P.~Shmerkin.
 Projections of self-similar and related fractals: a survey of recent developments, in
{\em Fractal Geometry and Stochastics V}, C. Bandt, K. Falconer and M. Z\"{a}hle (eds.), pp 3--25,
{\it  Progress in Probability.} {\bf 70}, Birkhauser, Cham, 2015.
 
\bibitem{Shm2}
P.~Shmerkin.
On the Hausdorff dimension of pinned distance sets,
{\it Israel J. Math.} {\bf 230} (2019), 949--972.

\bibitem{Shm3}
P.~Shmerkin.
A non-linear version of Bourgain's projection theorem,
{\it  J. Eur. Math. Soc.} {\bf 25} (2023), 4155--4204.

\bibitem{SS}
P.~Shmerkin and B. Solomyak.
Absolute continuity of self-similar measures, their projections and convolutions,
{\it Trans. Amer. Math. Soc.} {\bf 368} (2016), 5125--5151.
 
\bibitem{SSou}
P. Shmerkin and V. Suomala.
 Spatially independent martingales, intersections, and applications,
 {\it Mem. Amer. Math. Soc.} {\bf 251} (2018), no. 1195.
 
\bibitem{SV}
K.~Simon and L.~V\'{a}g\'{o}. 
 Projections of Mandelbrot percolation in higher dimensions, to appear
in {\em Fractals, wavelets, and their applications}, pp. 175-190, {\it Springer Proc. Math. Stat.} {\bf 92}, Springer, Cham, 2014. 
 
\bibitem{Sol}
B. Solomyak.
Measure and dimension for some fractal families,
{\it  Math. Proc. Cambridge Philos. Soc.} {\bf 124} (1998), 531--546.

\bibitem{Tri}
C. Tricot. 
Two definitions of fractional dimension,
{\it Math. Proc. Cambridge Philos. Soc.} {\bf 91} (1982), 57--74.
 
\bibitem{Xia} Y. Xiao. Packing dimension of the image of
 fractional Brownian motion, {\em Statist. Probab. Lett. } {\bf 333} (1997), 379--387.

\end{thebibliography}
\end{document}